\documentclass{article}

\RequirePackage{amsthm,amsmath,amsfonts,amssymb}
\RequirePackage[numbers]{natbib}
\RequirePackage[colorlinks,citecolor=blue,urlcolor=blue]{hyperref}
\RequirePackage{graphicx}
\RequirePackage{authblk}

\theoremstyle{plain}

\newtheorem{theorem}{Theorem}[section]

\theoremstyle{definition}
\newtheorem{definition}[theorem]{Definition}

\theoremstyle{remark}

\theoremstyle{plain}
\newtheorem{proposition}[theorem]{Proposition}
\theoremstyle{definition}
\newtheorem{remark}[theorem]{Remark}
\newcommand{\nc}{\newcommand}
\nc{\Prob}{{\mathbb{P}}}
\nc{\para}{{\mathbb{T}}}\nc{\E}{\,{\mathbb{E}}\,}
\nc{\Nset}{{\mathbb{N}}} \nc{\Rset}{{\mathbb{R}}}\nc{\Zset}{{\mathbb{Z}}}
\nc{\eqd}{{\stackrel{d}{=}}}
\nc{\Var}{{\mathrm{Var}}}
\newcommand{\kgs}[1]{{\color{black} #1}}

\title{Bivariate change point detection in movement direction and speed}

\author[1]{Solveig Plomer}
\author[2]{Theresa Ernst}
\author[2]{Philipp Gebhardt}
\author[2]{Enrico Schleiff}
\author[1]{Ralph Neininger}
\author[1]{Gaby Schneider}
\affil[1]{Institute of Mathematics,
Goethe University Frankfurt}
\affil[2]{Institute for Molecular Cell Biology of Plants,
Goethe University Frankfurt}
\begin{document}
\maketitle

\begin{abstract}
Biological movement patterns can sometimes be quasi linear with abrupt changes in direction and speed, as in plastids in root cells investigated here. For the analysis of such changes we propose a new stochastic model for movement along linear structures. 
Maximum likelihood estimators are provided, and due to  serial dependencies of increments, the classical MOSUM statistic is replaced by a moving kernel estimator. Convergence of the resulting difference process and strong consistency of the variance estimator are shown. We estimate the change points and propose a graphical technique to distinguish between change points in movement direction and speed.
\end{abstract}




\section{Introduction}
	
	\subsection{Motivation} \label{sec:motivation} 		
	
	The description of movement patterns can be important for the understanding of various biological processes on multiple scales. One main goal is to understand potential causes of changes between such movement patterns \citep{Gurarie2016,Patterson2017}. This can be important for controlling the spread of infectuous diseases or invasive species, to understand impacts of climate change, e.g., on the change of migration routes \citep{Nathan2008,McClintock2014,Leblond2016,Teitelbaum2021} or to minimize the impact of new technology on the natural habitat \citep{Russell2014,Russell2016}.
	
	On the macro scale, differences and changes between movement patterns of animals are investigated in order to describe specific behavior such as foraging strategies and predation, dispersal and migration, habitat use, social and territorial behavior, the coexistence of competitors or community interactions \citep[for an overview see, e.g.,][]{Byrne2009,Edelhoff2016,Gurarie2016,Bailey2018}. For example, highly explorative behavior with fast changes in movement direction and speed can alternate with resting periods with little changes or local exploration with minimal movement speed but highly diverse directions. Animals studied in this context are for example wolves \citep{Gurarie2016}, whales \citep{Gurarie2017CorrelatedVM,VacquiGarcia2018}, seals \citep{McClintock2014,Bailey2018} or butterflies \citep{Bailey2018}.
	
	On the micro scale, cells for example are compartmentalized into different cell organelles that form specialized reaction rooms in order to dissect  various metabolic processes. At the same time, most metabolic and signalling pathways are shared between different cell organelles. Consequently, the distribution of nutrients and signals within the cell needs to be coordinated, which requires in part organelle interactions \citep{Block2015,Shai2016,PericoSparks2018,Wang2023}. This is one, but not the only driving force for organelle movement within cells. One mode for this movement is cytoplasmic streaming \citep{Zawadzki1986} which can be observed in most of the eukaryotic cells  \citep[e.g.,][]{Ganguly2012,Alberti2015}. Aside from this rather arbitrary force a more directed mechanism of movement is based on the structures provided by the cytoskeleton \citep{PericoSparks2018}.
	
	So far, we only start to understand the principles of movement of the different cell organelles, but we are far from understanding the various regulations in different tissues and cell types. While new technological developments in the field of light sheet microscopy and super resolution microscopy allow the experimental determination of the cellular dynamics, fast procedures for the analysis of the large data sets are required. Especially a change point analysis for velocity and direction of movement would provide information, e.g., on involvement of cytosceleton components, influence of other organelles and the frequency of the switch between different modes of movement.\\
	
\begin{figure}
		\centering
		\includegraphics[width=1\textwidth]{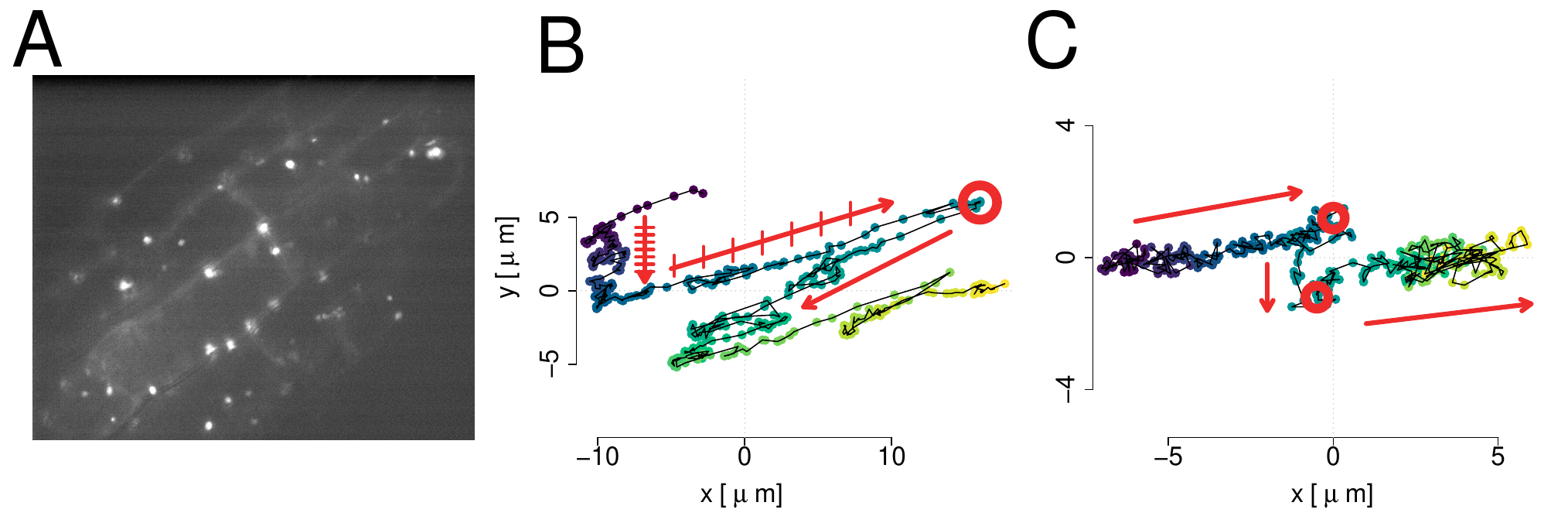}
		\caption{(A) Example of a microscopy image with fluorescently labelled plastids appearing as light spots. (B, C) Examples of two organelle tracks, no.~102 (B) and 434 (C) projected onto the first two principal components. Time is colour coded starting at dark blue. The phenomenon of changes in movement direction and step length is illustrated with red lines and arrows.}
		\label{fig:frame}
	\end{figure}

	In order to learn about the combination of and switch between different types of movement, we  consider here the movements of specific cell organelles in the root of the plant \textit{Arabidopsis thaliana} (see Figure \ref{fig:frame} A). These organelles are so-called plastides, which predominantly consist of leucoplasts with mainly nutrient storage functions and peroxisomes which take part in various reaction pathways. Details on the detection, tracking and post-processing of the two datasets can be found in Dryad \citep{dryad}.
	
	Examples of such tracks are illustrated in Figure \ref{fig:frame} B and C. Interestingly, while all tracks were recorded in three dimensions, over 90\% of the tracks showed more than 95\% of their movement variability in only two dimensions. We therefore focused on the two-dimensional projection of the movements, allowing comparability to approaches for animal movement pattern analysis. Our methods will, however, be also applicable to a higher number of dimensions.
	
	In the present data set, we observed interesting new features of the movement patterns (see Figure \ref{fig:frame} B and C): First, the movement of the cell organelles seemed to occur  along roughly linear structures (red arrows in panel B \kgs{and C}). Second, sharp changes in the movement direction were visible (red circles in panels B and C). Third, the speed, or the step length from one data point to the next, which are recorded in equidistant time intervals, showed a certain variability (red dashes in panel B). 
 
    It is therefore the goal of the present article to describe such  movement patterns in a mathematical framework and to derive a method that can statistically test for, detect and distinguish changes in the movement direction and speed in such tracks.

    \subsection{State of the art}\label{sect:stateoftheart}
	
	The offline change point analysis and detection of structural breaks is a major field of research \citep[see, e.g.,][]{Brodsky1993,Csorgo1997,AueHorvath2012,Jandhyala2013,Brodsky2017}. One considers abrupt changes in the mean or variance in $n$-dimensional real valued time series with independent observations \citep{Brodsky1993, Csorgo1997,Jandhyala2013}, cases with short or long term dependence \citep{AueHorvath2012, Dehlingetal2020, SchmidtFriedDehling2021,Padilla1,Padilla2} or cases with gradual changes \citep{VogtDette2015,Brodsky2017}. Change point detection methods range from likelihood ratio and Bayes-type methods \citep{Jandhyala2013} to CUSUM \citep{AueHorvath2012, Fryzlewicz2014} or MOSUM-based methods \citep{EichingerKirch2018, messerMFT, kirch2021moving}.
	
	However, only a small number of statistical methods is available for change point analysis in the field of directional statistics \citep[see e.g.,][]{Schnute1992,mardia2000directional,Jammalamadaka2001,Fitak2017}. When investigating change points in directional time series, classical change point analysis methods are often limited by  technical restrictions. For example, quantities corresponding to the standard deviation are rather limited, or the usual behavior of means or linear combinations of normally distributed random variables does not exist for the relevant directional distributions. Many approaches for change point detection in animal movement therefore focus on the partitioning of the observed time series, often lacking statistically rigorous tests for the null hypothesis of no changes in movement direction or speed \citep[see e.g.,][]{Byrne2009,Gurarie2009,Edelhoff2016,Gurarie2016,michelot2016}. In many cases, one intends to segregate  a small number of different states, where this number often needs to be specified in advance \citep{Lombard1986,michelot2016}. 
	
	Among the mathematical models applied for such analyses, one finds time discrete as well as time continuous models. In certain applications, e.g., with large sections of missing data or unequally spaced recording times, time continuous models such as Ornstein-Uhlenbeck processes or correlated velocity models can be necessary \citep[see][]{Johnson2008,Gurarie2017CorrelatedVM,Parton2017,Harris2013,McClintock2014,Fleming2017,Michelot2019}. In the present study, however, missing data were virtually not observed, all recording times were equally spaced, and the sampling rate was sufficiently small to capture important activity such as fast transport \citep{Tominaga2013}.
	
	Concerning time discrete models, a random walk (RW) is one of the most widely used families of models \citep[e.g.,][]{BaileyDissertation,Bailey2018,Bartumeus2008,Kareiva1983,McClintock2014,McClintock2012,michelot2016,Morales2004} \kgs{and applied in different variants. In animal studies, so-called correlated random walks \citep[CRWs, ][]{Kareiva1983,McClintock2012,michelot2016,Morales2004} are often applied in which the current movement  is assumed to depend on the current position and the previous movement direction. In such models, movement is described by choosing a random step size and turning angle at each movement step.   Turning angles around zero and high step sizes indicate strongly directed behavior as opposed to, e.g., explorative behavior with small step sizes. Such models do not directly fit the present data structure because they cannot maintain a global movement direction (Figure \ref{fig:Discussionmodelassumption} C). Also, changes in the movement direction during strongly directed behavior as we observe in the present data set (Figure \ref{fig:frame}) cannot be investigated in CRWs because the mean turning angle is always zero during directed behavior and does not depend on the absolute movement direction.} 

 \kgs{More suitable for the description of the observed movement along linear structures is therefore a so-called biased RW, which has been investigated, e.g., in the context of swimming micro-organisms \citep{Alt1980, Codling2010, Hill1997}. In such biased RWs, one assumes a mean absolute movement direction, which is distorted by some random component. Analyses using a Hidden Markov Model suggest that  assuming  absolute angles  in a biased RW may be more appropriate for the present data set than using turning angles as in CRWs \citep{Plomer2024b}. 
 
 The biased RW will therefore serve as a reference model in the present paper. As this biased RW is the only type of RW investigated here, we will refer to it simply as 'RW' in the following.} 
 
 The RW has technical advantages, e.g., it allows a relatively straightforward extension of existing change point methods. However, our findings indicate that a RW shows a higher variability in the movement direction and may thus move less strictly along a linear structure than can be observed in many organelle tracks \kgs{(Figure \ref{fig:Discussionmodelassumption} B). One possibility to deal with this issue would be to use more complex models, such as biased and correlated random walks, e.g., with attraction points. However, }such models can contain complex dependence structures and often need to be fitted using Bayesian methods or artificial neuronal networks \citep{Tracey2010,McClintock2012}. \kgs{Here, we follow an approach that proposes an alternative model in which the process sticks more closely to a linear structure than a RW.}

    \subsection{Aims and plan of the paper}	\label{sect:aims}
	
	In the present paper, we follow three main aims. First, we aim at developing a stochastic model for the description and statistical analysis of the recorded organelle tracks. This model will be presented in Section \ref{sec:model}. For the present recording setup, (i) the model can be time discrete, and it will even be constructed in such a way that the parametrization does not depend on the time resolution of the recording. Most importantly, the model should  be able (ii) to represent movements along roughly linear structures and (iii) to capture an arbitrary number of changes in (iv) the movement direction, or speed, or both. Due to these properties, the model will be called a Linear Walk (LW). As we also aim to capture changes  in the absolute movement direction, we propose to use a model with absolute angles (v).
	
	Second, we aim at developing a statistical test for the null hypothesis of no changes in direction or speed, which can also be used for the investigation of change points on multiple time scales. As we will show, changes in either direction or speed can be investigated individually, but only under the assumption that the respective other parameter is kept constant (Section \ref{sec:idea} and Supplement 1). If we allow for potential changes in both parameters, changes in either parameter will affect the analysis of the other. One would therefore need to know the changes in, e.g., the speed to analyse changes in the direction, and vice versa. We therefore propose a bivariate approach that investigates changes in these two parameters simultaneously.
 
    We then switch from the polar parameterization to a classical cartesian parameterization in which, for example, a classical MOSUM approach has been proposed for the bivariate analysis of  changes in  the mean and variance of a one-dimensional, independent and piecewise identically distributed time series by Messer \citep{messerBivariate}.  We illustrate its adaptation to our RW setting and show that it is less suitable in the LW setting due to a specific serial dependence structure in the increments (Section \ref{sec:MOSUMRW}).

    We  derive new estimators for the model parameters in the LW setting and use these to replace the classical MOSUM statistic with a moving kernel process (Section \ref{sec:MOSUMLW}). 	We then show that for a fixed window size under the null hypothesis, the new statistic scaled with the true variance is a Gaussian process whose distribution is independent from the model parameters (Section \ref{sec:MOSUMLW}). For increasing window size, we show convergence to a limit Gaussian process with a specific autocovariance structure (Section \ref{sect:konvergenz}). We also show strong consistency of the estimator $\hat{\sigma}^2$ for the variance $\sigma^2$ in the proposed model.

	Our third aim is to provide a method for change point estimation and classification of change points in direction or speed only, or in a combination of the two parameters. To that end, we adapt a change point detection algorithm that was proposed in the classical MOSUM setting, which allows testing for change points on multiple time scales and which has well studied asymptotic properties \cite{kirch2021moving,kirchreckruehm2024}. For classification of change points, we apply a graphical method called the \textit{leaf plot} in which changes in the two parameters typically result in specific graphical representations that represent different leaves. The methods are illustrated by application to tracks from the two sample data sets in Section \ref{sec:application}.
	
	\section{The Random Walk and the Linear Walk} \label{sec:model}
	
	\subsection{Notation and  model assumptions}\label{sect:Null}		
	
	Here we shortly review the Random Walk (RW) and propose a new stochastic model, the Linear Walk (LW). \kgs{When we refer to a LW or RW specifically, we will use the superscripts $LW$ and $RW$, respectively, to denote processes, variables or estimators. When the superscripts $LW$ and $RW$ are omitted, we refer to both, the LW and RW case.} \kgs{In both cases, we consider a bivariate} stochastic movement process $(X_i)$ with increments $Y_i:=X_i-X_{i-1},$  $i=1,2,\ldots$ The variables $Z_i, \tilde Z_i$ describe independent and bivariate standard normally distributed random variables. In both models, movement direction is denoted by $\vartheta \in [0,2\pi )$ and measured as an absolute angle. The step length is denoted by $r > 0$ and the variance of the error term is denoted by $\sigma^2>0$. 
	
\begin{figure}[htbp]
		\centering
		\begin{align}
			\begin{array}{cc}\includegraphics[width=0.8\textwidth]{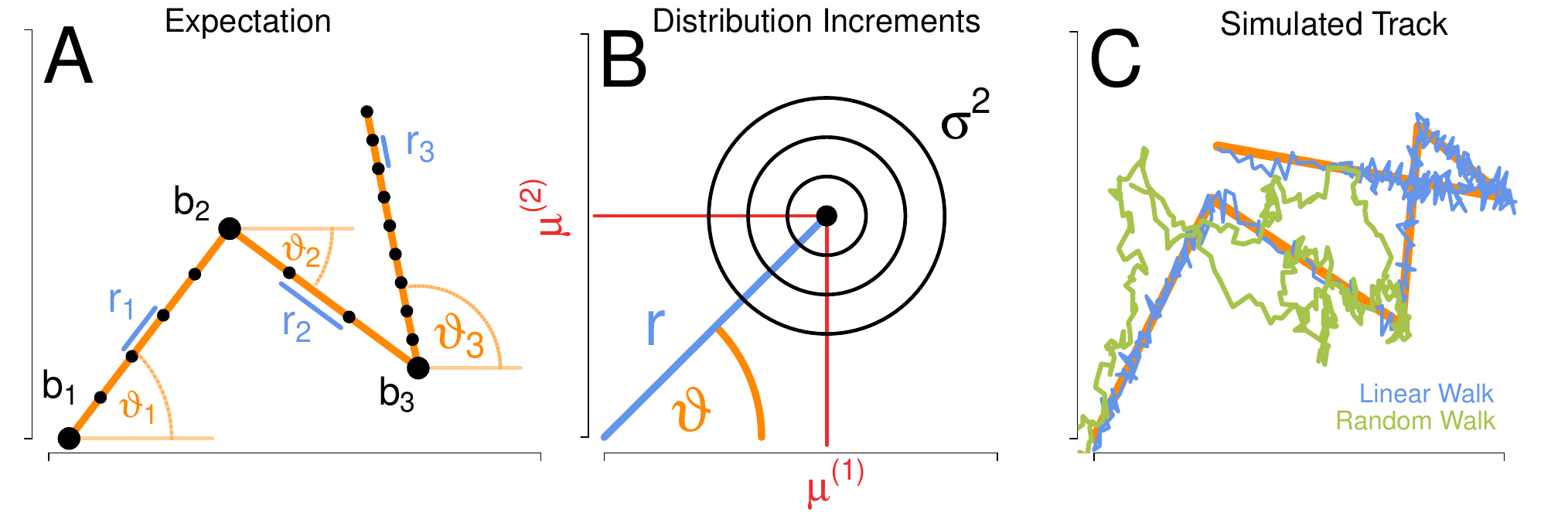}
			\end{array}\nonumber
		\end{align}
		\caption{(A) Visualization of the model parameters and the Expected Process (EP), (B) Polar and linear parametrization of the model, distribution of increments of the model, (C) visualization of the random displacement of the process by normally distributed error terms  in the LW (blue) and in the RW (green) around the EP (orange).}	
		\label{fig:model}
	\end{figure}

	In addition to the polar parametrization $(\vartheta,r)$, which allows a straightforward biological interpretation, we will also use the cartesian parameters $\mu:=(\mu^{\kgs{(1)}},\mu^{\kgs{(2)}})^t$ due to technical advantages, where (see {Figure \ref{fig:model} B})
	\begin{align}
		(\mu^{\kgs{(1)}},\mu^{\kgs{(2)}})^t= r(\cos \vartheta, \sin \vartheta)^t, \quad \vartheta := {\arctan}2(\mu^{\kgs{(2)}},\mu^{\kgs{(1)}})  \quad \text{and} \quad 
		r := \lVert \mu \rVert_2,  \label{def:ml_theta_r}
	\end{align}
	where $\lVert x \rVert_2$ denotes the euclidean norm and ${\arctan}2$ denotes a generalized inverse of the tangens function. Maximum likelihood estimators (MLEs) for $\vartheta$ and $r$ can therefore be derived from the MLEs of $\mu^{\kgs{(1)}}$ and $\mu^{\kgs{(2)}}$.

	\kgs{The RW and LW are formally defined in Definition \ref{def:RWLW}. Basically, the RW} assumes that all increments $Y_i\kgs{^{RW}}$ are independent and bivariate normally distributed, i.e., $Y_i\kgs{^{RW} = X^{RW}_i-X^{RW}_{i-1} = \mu + \sigma  \tilde Z_i, i\in \mathbb N}$. 
 
 \kgs{The parameter $\mu =(\mu^{{(1)}},\mu^{{(2)}})^t$, which is different from zero due to $r>0$, describes a biased movement in the direction given by $\vartheta = \text{arctan2}(\mu^{{(2)}},\mu^{{(1)}})$. For ease of notation, we refer to this {biased} RW simply as RW. The model assumptions and relations to other models in the context of movement analysis are discussed in more detail in Section \ref{sect:discmodelass}.} 
 
 \kgs{In comparison to the RW, we propose here a model called Linear Walk (LW) which has the same movement bias $\mu$ but sticks more closely to a straight line than the RW. Within a stationary section, we assume that an LW follows}
	\begin{align*}
		X^{LW}_{i} = b + i \mu + \sigma {Z_i},
	\end{align*}
\kgs{for $i \in \mathbb N$, with an arbitrary starting point  $b\in \mathbb R^2$.}

\kgs{Both an LW and RW can share the same process of expected values $e_i:=\mathbb E X_i^{RW}=\mathbb E X_i^{LW}$, differing only in the variability around this expectation. To see this, we define the expected process (EP) in Definition \ref{def:expectedwalk} and the corresponding LW and RW in Definition \ref{def:RWLW}}. In the full model with change points, we allow the EP $(e_i)\kgs{_{i\in \mathbb N}}$ to show changes in the movement direction $\vartheta$ and in the step length, $r$ (Fig.~\ref{fig:model} A).

\begin{definition}[\kgs{Expected Process, EP}]\label{def:expectedwalk}
\kgs{For a starting point $b\in \mathbb R^2$ and drift $\mu = (\mu^{(1)},\mu^{(2)})^t\neq 0 \in \mathbb R^2$, let $r$ and $\vartheta$ be as in \eqref{def:ml_theta_r}. Then the EP $(e_i)_{i\in\mathbb{N}}$ denoted by \textbf{EP}($\vartheta$,$r$,$b$,$\emptyset$) with an empty set of change points is given by  
$$e_i:= b + i \mu = b + i\cdot r (\cos\vartheta,\sin\vartheta)^t, \quad i\in \mathbb N.$$
An EP with a set of $k \in \mathbb N$ change points $C:=\{c_1,\ldots,c_{k}\}$ with $0=:c_0<c_1<\cdots<c_k<c_{k+1}:=T<\infty$ is defined as a composition of EPs without change points. To that end, define $(k+1)$ drifts $\mu_1,\ldots,\mu_{k+1}$, and let the sets of directions and step lengths be given by
$$\Theta = \{\vartheta_{1},\ldots,\vartheta_{k+1}\} \quad \text{and} \quad R = \{r_{1},\ldots,r_{k+1}\},$$ 
with
$\vartheta_j := {\arctan}2(\mu^{(2)}_j,\mu^{(1)}_j)$ and $r_j:=\lVert \mu_j \rVert_2$, $j=1,\ldots,k+1$, and define}  $k\kgs{+1}$ processes \textbf{EP}($\vartheta_j$,$r_j$,$b_j$,$\emptyset$) that describe the piecewise \kgs{linear} EPs in the sections between the change points. The \kgs{set of starting points $B=\{b_1,\ldots,b_{k+1}\}$}  are defined such as to connect the \kgs{EPs}, i.e., $b_1\in\mathbb{R}^2$ and 
$$b_{j}=(c_{j-1}-c_{j-2})\cdot r_{j-1}\cdot (\cos(\vartheta_{j-1}),\sin(\vartheta_{j-1}))^t+b_{j-1}, \quad j\in\{2,\ldots,k\kgs{+1}\}.$$ 
We assume that for successive EPs, at least one of the parameters $\vartheta$ or $r$ changes. 
\kgs{Then $(e_{i}^{(j)})_{i\in\mathbb{N}}=$ \textbf{EP}($\vartheta_{j}$,$r_{j}$,$b_{j}$,$\emptyset$) are given by } $e_{i}^{(j)}\kgs{:=b_j+i\mu_j}$ \kgs{, and $(e_i^C)_{i=1,\ldots,T}=$} \textbf{EP}($\kgs{\Theta,R,B}$,C) with change point set $C$ is defined as
		\begin{align}
			e_i^C:=e^{(j\kgs{+1})}_{i-c_j},\quad j=0,\ldots,k, \quad i=c_j+1,\ldots,c_{j+1}. \label{eq:expectedprocesscomposite}
		\end{align}
	\end{definition}

	\kgs{\begin{definition}[RW and LW]\label{def:RWLW}
Let  $(e_{i}^{C})_{i\in\mathbb{N}}=$ \textbf{EP}($\Theta$,$R$,$B$,$C$) denote an expected process with parameter sets $\Theta,R,B$ and change point set $C$ as in Definition \ref{def:expectedwalk}. Let $\sigma^2>0$ and $Z_i, \tilde Z_i, i=1,2,\ldots,$ be independent and bivariate standard normally distributed random variables. Then the corresponding RW and LW, \textbf{RW}($\Theta,R,B,\sigma^2$,C) and \textbf{LW}($\Theta,R,B,\sigma^2$,C), are given, respectively, by 
\begin{align}
X_i^{RW}&:=e_i^C + \sigma \sum_{j=1}^i \tilde Z_j \quad\quad \text{and}\label{eq:RW}\\
X_i^{LW}&:=e_i^C + \sigma Z_i, \quad \text{for }i\in\mathbb N.  \label{eq:LWbasic}
\end{align}
		\end{definition}}
  
    \subsection{Discussion of the model assumptions}\label{sect:discmodelass}
    The RW and LW models capture the two features of the absolute movement direction and step length with two biologically interpretable parameters. Changes in these parameters can be implemented separately, which enables the description of  different modes of movement. For example, a high step length and  constant movement direction could indicate transport along intracellular filaments, while lower step length could indicate cytoplasmic streaming. Note that we use normally distributed errors because these have technical advantages over generic directional distributions.
	
\begin{figure}[htbp]
		\centering
		\begin{align}
			\begin{array}{cc}\includegraphics[width=0.95\textwidth]{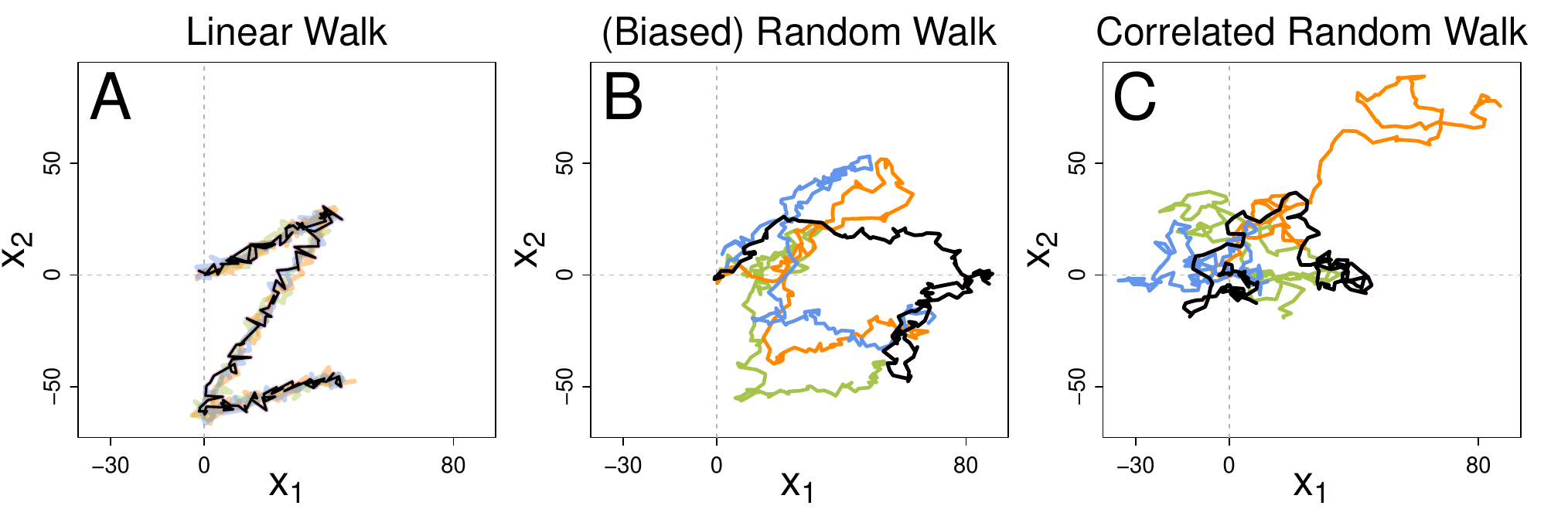}
			\end{array}\nonumber
		\end{align}
		\caption{\kgs{Four simulated tracks of length $T=150$ within the three movement models LW (A), (Biased) RW (B) and CRW (C). Parameters for all models were $\vartheta=(35,-115,20)/180\cdot\pi$, $r=(1,2,1)$, $\sigma=2$, $b=(0,0)$ and $C=\{50,100\}$.  LW and RW were simulated as in Definition \ref{def:RWLW}. For the CRW, for comparability, we used the same parameters, change points and assumptions on the error terms, but the turning angles were set to zero in sections of constant movement direction, while at each change point, one fixed turning angle was chosen.}}
		\label{fig:Discussionmodelassumption}
	\end{figure}

	Qualitatively, the LW (Fig.~\ref{fig:model} C, blue, \kgs{and Fig.~\ref{fig:Discussionmodelassumption} A}) sticks more closely to the EP (\kgs{Fig.~\ref{fig:model} C, orange, and Fig.~\ref{fig:Discussionmodelassumption} A}) and thus captures the linear movement of the plastid tracks (Fig.~\ref{fig:frame}). If the EP is interpreted as a filament along which the organelle is transported, the random displacements $Z_i$ of the points assumed in the LW could be interpreted as a fluctuation of the organelle around its adhesion point or as a measurement error. In an RW (Fig.~\ref{fig:model} C, green\kgs{, and Fig.~\ref{fig:Discussionmodelassumption} B}), the linear structure hidden in the EP is much less visible due to a cumulation of error terms $\tilde Z_i$. \kgs{In a CRW (Fig.~\ref{fig:Discussionmodelassumption} C), linear movement can only hardly be observed for the present parameter values.}

	The RW and LW also differ with respect to the dependence structure in their increments. While the increments in the RW are independent and $\mathcal N(\mu,\sigma^2 I)$-distributed by definition, the increments in the LW are given by 
 \begin{equation*}
        Y^{LW}_{i}= \mu + \sigma (Z_{i}-Z_{i-1}).
\end{equation*} 
Thus, they exhibit a serial dependence of order one because successive increments $Y_i$ and $Y_{i+1}$ share the error term $Z_i$. For each dimension $d=1,2,$ we thus observe
	\begin{align}
		\text{Cor}(Y^{LW\kgs{,(d)}}_i,Y^{LW\kgs{,(d)}}_{i+1}) = -1/2, \quad  \text{Cor}(Y^{LW\kgs{(d)}}_{i},Y^{LW\kgs{,(d)}}_{i+k}) = 0 \quad \text{for} \quad k\ge 1. \label{eq:covarianceLinearWalk}
	\end{align}
	In Section \ref{sect:MLE} we will discuss the behavior of the classical mean of increments, which shows crucial differences between the RW and the LW.

    \subsection{Parameter estimation}\label{sect:MLE}
	
	In order to derive the maximum likelihood estimates (MLEs) of the model parameters, both for the RW and the LW, we derive the estimators for the parametrization $(\mu^{(1)},\mu^{(2)})^t$ and then obtain the estimators for the direction $\vartheta$ and step length $r$ using equation \eqref{def:ml_theta_r}.  	
	
	In the RW, the increments $Y^{\kgs{RW}}_i$ are independent and normally distributed, so that the MLEs can be  derived classically as the mean and empirical variance. Thus, as estimates in a window of size $h\in\mathbb{N}, h>2$ starting at $i$, we  obtain in the RW
	\begin{align*}
		\hat{\mu}^{RW}(i,h)&=\frac{1}{h}\sum_{j=1}^{h}Y^{\kgs{RW}}_{i+j} \quad \text{and} \\
  {(\hat{\sigma}^{RW})}^2(i,h)&=
  \kgs{\sum_{j=1}^h\|Y_{i+j}^{RW}-\hat\mu^{RW}(i,h) \|^2}.
	\end{align*}
	
    Interestingly, for the LW, the classical mean of increments behaves quite differently than in the RW, and it does no longer represent the MLE for $\mu$. In the RW, the variance of $\hat{\mu}^{\kgs{RW}}_h\kgs{(i,h)}$ falls like $h^{-1}$ due to the independence of increments, i.e., for every dimension $d=1,2,$
    \begin{equation*}
\Var(\hat{\mu}^{RW,(d)}\kgs{(i,h)})=h^{-1}{\Var (Y^{\kgs{RW,(d)}}_{i})}= h^{-1}{\sigma^2}.
\end{equation*}
In contrast, in the LW, the classical mean reduces to a telescope sum, i.e.,  $\frac{1}{h}\sum_{j=1}^hY^{\kgs{LW}}_{j} = \mu + \frac{\sigma}{h} (Z_{h}-Z_{1}),$	implying that its variance,
${2\sigma^2}h^{-2}$,
 falls like $h^{-2}$ within the LW. As a consequence, the classical MOSUM statistic will be replaced by a kernel estimator in Section \ref{sec:MOSUMLW}, replacing the classical mean by the maximum likelihood estimator (MLE) of $\mu$ as given in Proposition \ref{theorem:ml_mu_positions}.	
	
	\begin{proposition}\label{theorem:ml_mu_positions}
		Let \textbf{LW}($\vartheta,r,b$,$\sigma$,$\emptyset$) be a uni-directional LW as in \kgs{Definition} \ref{def:RWLW}, and let \kgs{for $i\in\mathbb N$},  $X^{\kgs{LW}}_{i+1}, X^{\kgs{LW}}_{i+2}, \ldots, X^{\kgs{LW}}_{i+h}$ be a sequence of positions from this model. Then the MLEs for $\mu, b$ and $\sigma^2$ are given, respectively, by 
		\begin{align}
			\hat{\mu}^{\kgs{LW}}(i,h)&=\frac{6}{h^3-h}\sum_{j=1}^{\lfloor\frac{h}{2}\rfloor}\left(h-(2j-1)\right)\left(X^{\kgs{LW}}_{i+j}-X^{\kgs{LW}}_{i+h-j+1}\right), \label{eq:MLEmuLM}\\
			\hat{b}^{\kgs{LW}}(i,h)  &= \bar X\kgs{_i}^{\kgs{LW}}- (\kgs{i}+(h+1)/2)\hat\mu^{\kgs{LW}}\kgs{(i,h)} \quad \text{and} \label{eq:MLEbLM}\\
			(\hat{\sigma}^{\kgs{LW}})^{2}(i,h) &=\frac{1}{2h}\sum_{j=1}^{h} \| X^{\kgs{LW}}_{i+j}-(i+j)\hat\mu^{\kgs{LW}}\kgs{(i,h)}-\hat b^{\kgs{LW}}\kgs{(i,h)}\|^2,
			\label{eq:MLEsigmaLM}\end{align}
	where $\bar{X}\kgs{_i}^{\kgs{LW}}:=h^{-1}\sum_{j=1}^h{X^{\kgs{LW}}_{i+j}}$ denotes the mean of the observed points in the respective window.
	\end{proposition}

	The proof can be found in Appendix \ref{proof:ml_mu_sigma_positions_t}. The MLEs of $\mu$ and $b$ are unbiased. In the estimator of $\sigma^2$ we replace the numerator $2h$ by $2h-4$ to obtain an unbiased version of the estimator (see Supplement 2). The MLEs of the direction\kgs{, $\hat \vartheta(i,h)$,} and step length\kgs{, $\hat r(i,h)$,} are given as functions of $\hat \mu(i,h)$ according to equation \eqref{def:ml_theta_r}.

	In the LW, the MLEs of $\mu$ have an interesting geometric interpretation in terms of the increments $(X_{i+j}^{\kgs{LW}}-X_{i+h-j+1}^{\kgs{LW}})$ between the $j$th and the $j$-last observation in equation \eqref{eq:MLEmuLM}. The largest increment from the first to the last position within the window, $X_{i+1}^{\kgs{LW}}-X_{i+h}^{\kgs{LW}}$, gets the largest weight $(h-1)$, and the weight decreases linearly with the temporal distance between positions. \kgs{As a potential intuitive explanation of this result, all these increments have the same variance $2\sigma^2$, but the expectation of the increments in equation \eqref{eq:MLEmuLM} and thus, the relative precision, is largest for $j=1$, i.e., for the largest increment from the first to the last observation, and decreases linearly with $j$.} Also the parameter $b$ is estimated intuitively by starting with the center of mass ${\bar X}$ and subtracting an appropriate number of steps of size $\hat \mu\kgs{^{LW}(i,h)}$.
	
	In the following, we will develop methods for change point detection within the two models, with a focus on the newly proposed LW.

 \section{Statistical test and change point detection}\label{sec:cp}
 
	In this section we derive a statistical test for the null hypothesis of no change in the movement direction or speed \kgs{both for the RW and the LW}. Section \ref{sec:idea} will present the basic idea and explain the limits of two seemingly straightforward approaches, namely the univariate analysis and the classical MOSUM approach.

	\subsection{The idea and limitations of seemingly natural approaches}\label{sec:idea}
\begin{figure}[htbp]
		\centering
		\includegraphics[width=1\textwidth]{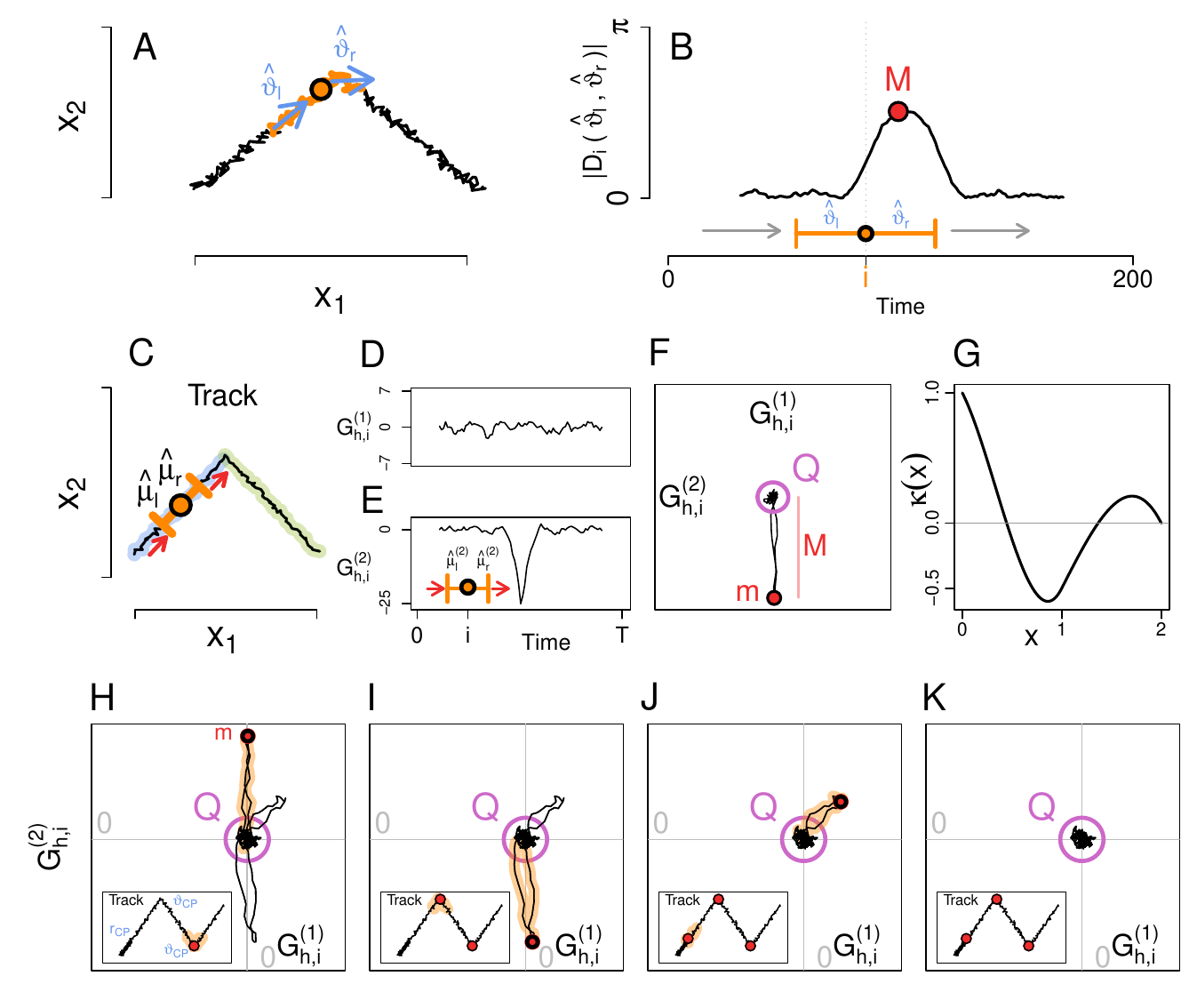}		
		\caption{(A) Realization of a LW track with one change in the direction. Double window (orange line) around time $i$ (black circle) with $h = 30$. Blue arrows indicate estimated direction in left and right window. (B) Resulting process of absolute differences  $|D_i|$ and its maximum (red point). Double window (orange) around time $i$ and its moving direction (grey arrows). (C) A LW track (black line). Track before change point in direction (blue) and after (green). In a double window of size $2h$ around time $i$ (orange), estimates of the expectation in the left and right window\kgs{, $\hat\mu_{\ell}^{LW}:=\hat\mu^{LW}(i-h,h)$ and $\hat\mu_{r}^{LW}:=\hat\mu^{LW}(i,h)$,} are obtained. Red arrows indicate moving window. (D, E) \kgs{The two components of the p}rocess of scaled differences $(G_{h,i})_{i\in\mathbb N}$ and (F) $(G_{h,i})_{i\in\mathbb N}$ as a two-dimensional process. Point of maximal deviation from origin m (red dot), length of maximal deviation M (red line) and boundary of rejection circle with radius $Q$ (purple). (G) Normalized autocovariance function $\kappa$. (H--K) Change point detection algorithm. (H) Process $(G_{h,i})$ derived from a LW with three change points, rejection threshold (purple), point of maximum deviation, $m$ (red), and $2h$-neighbourhood (orange). The LW track is shown in the bottom left together with the position of $m$ (red). (I,J) The process $(G_{h,i})$, where the $2h$-neighbourhood around $m$ (I) and the second maximal deviation (J) has been deleted, new maximal deviations marked in red. (K) $(G_{h,i})$ remains within the rejection borders, CPD is finished.}
		\label{fig:biv1}
	\end{figure}

 	The main idea to statistically test for and detect changes in the movement direction and speed is to apply a moving window technique. By moving a double window of size $2h$ across \kgs{a track of length $T$}, one successively estimates the parameter of interest, for example the movement direction $\vartheta$, as $\hat{\vartheta}_{\ell}\kgs{(i):=\hat\vartheta(i-h,h)}$ and $\hat{\vartheta}_{r}\kgs{(i):=\hat\vartheta(i,h)}$ in the left and right half of this window (Figure \ref{fig:biv1} A, blue directions and orange window). One then calculates the process of differences, which would, in the case of $\vartheta$ denoting the direction, be defined as (cmp.~Figure \ref{fig:biv1} B), \kgs{for $i=h,\ldots,T-h$,}
	\begin{equation}
D_i:=D_i(\hat{\vartheta}_{\ell}\kgs{(i)},\hat{\vartheta}_{r}\kgs{(i)}):=
		\text{atan2}(\sin(\hat{\vartheta}_{r}\kgs{(i)}-\hat{\vartheta}_{\ell}\kgs{(i)}),\cos(\hat{\vartheta}_{r}\kgs{(i)}-\hat{\vartheta}_{\ell}\kgs{(i)})), \label{eq:mosum}
	\end{equation}
	which gives the smaller angle enclosed by the two directions. If all parameters are constant, the process  $(D_i)_{i\kgs{\in\{h,\ldots,T-h\}}}$ should fluctuate around zero, while it will show systematic deviations from zero in the neighbourhood of change points. Therefore, the maximum of the resulting absolute difference process,
	\begin{equation}
		M:=\max_i |D_i|,
	\end{equation}
	can serve as a test statistic for the null hypothesis of no change. A suitable quantile $Q$ of the distribution of $M$ under the null hypothesis can be used to define the rejection threshold of the statistical test.
	
    \subsubsection{The univariate approach}
	It would be intriguing to apply the described approach directly to univariate parameters such as the movement direction. However, one should note that this univariate approach has a number of problems and restrictions. First, if the difference process is defined as above for directions, the distribution of $M$ is not easily accessible formally. \kgs{Both for the LW and RW, t}he estimator $\hat{\vartheta}\kgs{(i,h)}$ follows a Projected Normal (PN) distribution because it is derived from the normally distributed estimator of the expectation $\hat{\mu}\kgs{(i,h)}$ by projection onto the unit circle (see e.g.~\cite{mardia2000directional}). While statements about the projected normal distribution are possible in simple scenarios \citep{WangGelfand2013,Mastrantonio2015,NunezAntonio2011}, we are not aware of general properties of either the distribution of the difference of two PN variables, nor process theory with PN distributed components. Therefore, the rejection threshold $Q$ needs to be derived in simulations based on the parameter estimates of all other process parameters. As a second restriction, this automatically requires the other parameters to be constant or their change points to be known beforehand (cmp.~\cite{Albert2017}). This is because changes in, for example, $r$, will lead to biased estimation of $r$ and thus, a biased threshold $Q$, such that the significance level cannot be kept. Finally, even if all model parameters are assumed constant, the step length $r$ should not be too small as compared to $\sigma$ because otherwise, estimation of $r$ will be biased. For more details on this approach see Supplement 1. It will therefore be important to consider both parameters simultaneously in a bivariate setting.

    \subsubsection{The classical MOSUM}\label{sec:MOSUMRW}
	
    As another idea, one might consider using a classical MOSUM approach, which has been proposed in a bivariate framework by Messer \cite{messerBivariate}.  While in the RW setting, this approach is straightforward and will be shortly recalled and adapted, it requires careful consideration and extension in the LW setting.

    In this approach, we make use of the parametrization $\mu=(\mu^{(1)}, \mu^{(2)})$ instead of the polar parametrization  $(\vartheta, r)$. Specifically, we apply a double window of size $2h$ centered around time $i$ to estimate the difference $\hat{\mu}^{\kgs{RW}}\kgs{(i,h)}-\hat{\mu}^{\kgs{RW}}\kgs{(i-h,h)}$. This results in a two dimensional process of differences, which we scale with its estimated variance, such that the derived test statistic will be independent of the model parameters even if changes in both parameters can occur.

    In short, we consider an \textbf{RW}(\kgs{$\Theta$,$R$,$B$},$\sigma$,$C$) as in Definition \ref{def:expectedwalk}. For a window of size $h\in\mathbb{N},h>2$ and $i=1,2,\ldots$, the difference of means is scaled with its variance, i.e.,
	\begin{align}
		G_{h,i}^{RW} &:= \frac{\sqrt{h}\left(\hat{\mu}^{\kgs{RW}}(i,h)-\hat{\mu}^{\kgs{RW}}(i-h,h)\right)
		}{\sqrt{{(\hat{\sigma}^{\kgs{RW}})^2(i-h,h)+(\hat{\sigma}^{\kgs{RW}})^2(i,h)}}}. \label{eq:ghtRW}
	\end{align}
    Note that for the RW, we find the classical mean in the numerator. Then, for an observed process of length $T\in\mathbb{N}$, we consider the two dimensional process of differences 
    $$G^{RW}:=\left(G_{h,i}^{RW}\right)_{i\in \{h,h+1,\ldots,T-h\}}.$$ 
    The process $G^{RW}$ will typically fluctuate around zero but show systematic deviations in the neighborhood of a change point. Thus, the maximal deviation from the origin,
	\begin{align}
		M_h:= \max_{i} \|G_{h,i}^{RW}\|^2,
	\end{align}
    is used as test statistic for the null hypothesis of no change in bivariate expectation, i.e., no change in  direction or step length. 

    In order to derive a rejection threshold for the statistical test, one can use straightforward modifications of the proof presented in \cite{messerBivariate}. For asymptotics a growing number of data points within the windows is required. To that end one switches from the time discrete setting with $i,h \in \mathbb N$ to a time continuous setting, where $i$ and $h$ are replaced by $nt$ and $n\eta$ respectively, where $\eta>0$ and $t\in[\eta,T-\eta]$ are real, $T>2\eta$, and $n\to\infty$ is considered. \kgs{All parameter estimates can be adapted to the time continuous setting, e.g., defining $\hat\mu^{RW}(nt,n\eta):=\hat\mu^{RW}(\lfloor nt \rfloor, n\eta)$, and analogously for the other terms. For convenience we assume $\eta\in\mathbb N$, which is not a restriction in practice because $\eta$ can be assumed a multiple of the recorded time resolution. All results of the present paper also hold, with slight adaptations, for  real $\eta>0$.}

    The process $(G_{n\eta,nt}^{RW})_{t\kgs{\in[\eta,T-\eta]}}$ can then be rewritten as a functional of sums of the form $\kgs{S^{(n)}(t):=n^{-1/2}\sum_{j=1}^{\lfloor nt\rfloor}\tilde{Z}_{j}^{(d)}}$. According to the functional limit theorem $\left(S^{(n)}(t)\right) \to \left(W_t\right)$ in distribution, which is, in various forms, a key component in MOSUM and CUSUM frameworks \citep{messerMFT,messerBivariate,AueHorvath2012,kirch2021moving} and can be used to show convergence. Then the process $(G_{n\eta,nt}^{RW})_{\kgs{t\in[\eta,T-\eta]}}$  convergences in distribution to a limit process $(\mathcal{L}_{\eta,t}^{RW})_{\kgs{t\in[\eta,T-\eta]}}$, which is a functional of a planar Brownian motion. That means, as $n \to \infty$,
	\begin{align}
		\left(G_{n\eta,nt}^{RW}\right) \xrightarrow[\quad\quad]{d} \left(\mathcal{L}_{\eta,t}^{RW}\right) \text{ in } \left(\mathcal{D}_{\mathbb{R}^2}[\eta,T-\eta],d_{SK}\right),\label{lhtRW}
	\end{align}
	where the two dimensional limit process  $(\mathcal{L}_{\eta,t}^{RW})$ does not depend on the model parameters. Here, $\left(\mathcal{D}_{\mathbb{R}^2}[\eta,T-\eta],d_{SK}\right)$ denotes the space of $\mathbb{R}^2$-valued c\`{a}dl\`{a}g-functions on $[\eta,T-\eta]$ equipped with the Skorokhod topology. The \kgs{vectors} $\mathcal{L}_{\eta,t}^{RW}$ are given by
	\begin{align}
		\mathcal{L}_{\eta,t}^{RW} = (2\eta)^{-1/2}((W_{t+\eta}-W_t)-(W_t-W_{t-\eta})),
	\end{align} 
	where \kgs{$(W_t)_{t\in[0,T]}:=(W^{(1)}_t,W^{(2)}_t)_{t\in[0,T]}$} denotes a planar Brownian motion. Using this limit process, a rejection threshold can be obtained via simulation cmp.~\cite{messerBivariate}. \kgs{For this we use $\sum_{j=1}^{\lfloor t\rfloor}\tilde Z_j$ as a discrete approximation of $W_t$ such that the process with elements\kgs{, for $t\in[\eta,T-\eta]$},}
    \begin{align*}
         \kgs{\Gamma^{RW}_{n\eta,nt}} & \kgs{:= (2 n\eta)^{-1/2}{\sum_{j=1}^{n\eta} \left(\tilde Z_{\lfloor nt\rfloor+j}-\tilde Z_{\lfloor nt\rfloor-n\eta+j}\right)}}
    \end{align*}
    \kgs{is an approximation of $(\mathcal{L}_{\eta,t}^{RW})$. In case of a change point, the process diverges, such that consistency of the test introduced in Section \ref{sect:practice} can be established, cf. \citep{kirch2021moving}.} 

 \kgs{\subsubsection{Limitations of the classical MOSUM in the LW}}
    \kgs{Whereas} the classical MOSUM is \kgs{a suitable tool} in the RW setting\kgs{, it would not be consistent in the LW setting, i.e., the testpower would not increase with $n$}. \kgs{This can be seen noting that the  variance of the classical mean decreases with $(n\eta)^{-2}$ instead of $(n\eta)^{-1}$ in the LW (Section \ref{sect:MLE}). Therefore, the MOSUM statistic would need to be scaled with $n\eta$ instead of $\sqrt{n\eta}$, yielding a MOSUM statistic of the form}
	\begin{align}
		&\kgs{\frac{n\eta}{\sqrt{2\sigma^2}}\left({\frac{1}{n\eta}\sum_{j=1}^{n\eta}Y^{LW}_{\lfloor nt\rfloor+j}-\frac{1}{n\eta}\sum_{j=1}^{n\eta}Y^{LW}_{\lfloor nt\rfloor -n\eta+j}}\right)}\label{eq16}\\
        &\quad\quad\quad \kgs{=} \kgs{\frac{\sigma}{\sqrt{2\sigma^2}}\sum_{j=1}^{n\eta}\left((Z_{\lfloor nt\rfloor+j}-Z_{\lfloor nt\rfloor+j-1}) - (Z_{\lfloor nt\rfloor -n\eta+j}-Z_{\lfloor nt\rfloor -n\eta+j-1}) \right),}\nonumber
    \end{align} 
    \kgs{where we used $Y^{LW}_{i}= \mu + \sigma (Z_{i}-Z_{i-1})$. Summing up the two inner brackets separately over $j$ yields telescopic sums, where the first sum reduces to $Z_{\lfloor nt\rfloor+n\eta}-Z_{\lfloor nt\rfloor}$ and the second to $Z_{\lfloor nt\rfloor}  - Z_{\lfloor nt\rfloor -n\eta}$. This results in \eqref{eq16} being equal to}
    \begin{align}
        &\frac{1}{\kgs{\sqrt 2}}\left(Z_{\kgs{\lfloor nt \rfloor+n\eta}}-2Z_{\lfloor nt\rfloor}+Z_{\kgs{\lfloor nt \rfloor- n\eta}}\right)\kgs{,}
    \end{align} 
	\kgs{which} has the distribution $N(0,\kgs{3 I_2})$ for every $n$, $\eta$ and $t$\kgs{, where $I_2$ denotes the identity matrix in $\mathbb R^2$}. As \kgs{the statistic no longer depends on $n$ or $\eta$}, the testpower cannot be increased by increasing the window size, such that small changes will always remain unlikely to be detected.	
 
    Therefore, we propose to replace the classical MOSUM by a moving kernel estimator that makes use of the MLE of $\mu$ derived within the LW. \kgs{For a note on the consistency of the resulting statistical test see Section \ref{sect:practice}.}\\

    \subsection{The \kgs{kernel} approach for the LW}\label{sec:MOSUMLW}
			
    For the LW setting, we propose to substitute the classical mean by its MLE of $\mu$ in  the derivation of the process $G$. In the time discrete setting, as an analogue of equation \eqref{eq:ghtRW} we thus consider the process
    \begin{align}
        G_{h,i}^{LW} := \frac{\hat{\mu}^{LW}(i,h)-\hat{\mu}^{LW}(i-h,h)}{\sqrt{\frac{12}{h^3-h}((\hat\sigma^\kgs{LW})^2(i-h,h)+(\hat\sigma^\kgs{LW})^2(i,h))}}, \quad \kgs{ i\in\mathbb N}.\label{eq:ghtLW}
    \end{align}
    The denominator of equation \eqref{eq:ghtLW} estimates the standard deviation of the numerator, because it holds
    \begin{align*}
        \mathbb{V}ar\left(\hat{\mu}^{\kgs{LW},(d)}(i,h)\right) &= \mathbb{V}ar\left( \frac{6}{h^3-h}\sum_{j=1}^{h}(2j-h-1)X_{i+j}^{\kgs{LW},(d)}\right)\\
        &=\frac{36}{(h^3-h)^2} \sum_{j=1}^n \left(2j-h-1\right)^2\mathbb{V}ar\left(X_1^{\kgs{LW},(d)}\right)\\
        &= \frac{12}{h^3-h}\sigma^2,
    \end{align*}
	where we use the independence of $X_i^{\kgs{LW},(d)}$ and $X_k^{\kgs{LW},(d)}$ for $i\neq k$.

    We study the process $(G_{h,i}^{LW})_\kgs{i\in\mathbb N}$ to construct a statistical test. For asymptotics we replace $i$ and $h$ by $nt$ and $n\eta$ respectively, where $\eta>0$ and $t\in[\eta,T-\eta]$ are real and $T>2\eta$ is a finite time horizon.  All parameter estimates can be adapted to the time continuous setting analogously to Section \ref{sec:MOSUMRW}, e.g., defining $\hat\mu^\kgs{LW}(nt,n\eta):=\hat\mu^\kgs{LW}(\lfloor nt \rfloor,n\eta)$, and analogously for the other terms\kgs{, where we again assume $\eta\in\mathbb N$ for convenience}.

    Firstly, in Proposition \ref{lemma:consistencyhatsigmasquared}, we obtain that the estimate of $\sigma^2$ from equation \eqref{eq:MLEsigmaLM} is pointwise strongly consistent. Secondly, we investigate the limit behavior of a time continuous version of $(G_{h,i}^{LW})_\kgs{i\in\mathbb N}$ in Section \ref{sect:konvergenz} with respect to a functional limit law, see Proposition \ref{proposition:convergenceupsilonn}.

    \begin{proposition}[Pointwise Strong Consistency of MLEs in the LW]\label{lemma:consistencyhatsigmasquared}
    Let \textbf{LW}({$\vartheta,r,b$},$\sigma$,$\emptyset$) be a uni-directional LW as in equation (\ref{eq:LWbasic}) and let \\$X^{\kgs{LW}}_{1},X^{\kgs{LW}}_{2}, \ldots$ be a sequence of positions from this model. Let $\eta,t>0$. Then the MLEs $\hat{\mu}^\kgs{LW}\kgs{(nt,n\eta)}$, $\hat{b}^\kgs{LW}\kgs{(nt,n\eta)}$ and $(\hat{\sigma}^\kgs{LW})^2\kgs{(nt,n\eta)}$ as stated in Proposition \ref{theorem:ml_mu_positions} are pointwise strongly consistent, i.e., as $n\rightarrow\infty$ almost surely
		\begin{align}
			\hat{\mu}^\kgs{LW}(nt,n\eta) \to \mu,\quad\quad
			\hat{b}^\kgs{LW}(nt,n\eta) \to b,  \quad\quad
			(\hat{\sigma}^\kgs{LW})^2(nt,n\eta) \to \sigma^2.\label{eq:strongconsistency}
		\end{align}
	\end{proposition}
	For the proof, the Strong Law of Large Numbers (SLLN) is adapted to a setting in which each summand is multiplied by a weight dependent on $n$ and the index of the summand. One version of the SLLN for weighted sums and the proof of the strong consistency for $\hat{\mu}^\kgs{LW}\kgs{(nt,n\eta)}$ can be found in Appendix \ref{appendix:SLLNnew}. The proof of the strong consistency of $(\hat{\sigma}^\kgs{LW})^2\kgs{(nt,n\eta)}$ follows the same structure but involves more tedious calculations. Together with the second required SLLN and the proof for the strong consistency of $\hat b^\kgs{LW}\kgs{(nt,n\eta)}$, it can be  found in Supplement 3.\\

    Now we define a time continuous version of \eqref{eq:ghtLW} to investigate process convergence. 
    \begin{definition}[Process of Differences $\Gamma_{h,t}^{\kgs{LW}}$]\label{def:ght}
		Let $\left(X_{i}^{\kgs{LW}}\right)_{i\in \mathbb N_0}$ be a sequence of positions from an \textbf{LW}(\kgs{$\Theta,R,B$},$\sigma$,$C$), and let $t\in\mathbb{R}$ and $\eta\in\mathbb{R}, \eta>0$. Let $\hat{\mu}^\kgs{LW}(i,h)$ be the MLEs of $\mu$ in the LW setting according to Proposition \ref{theorem:ml_mu_positions}. 
        Then, for $n\in\mathbb{N}$,  let
		\begin{align}
			\Gamma^{\kgs{LW}}_{n\eta,nt} &:= \frac{\hat{\mu}^\kgs{LW}(nt,n\eta)-\hat{\mu}^\kgs{LW}(n(t-\eta), n\eta)}{\sqrt{\frac{24}{\kgs{(n\eta)}^3-\kgs{n\eta}}\sigma^2}}\label{eq:GammaLW}.
		\end{align} 
    \end{definition}
    Note that the two components of the process are stochastically independent by construction.
    Under the null hypothesis, cf.~(\ref{eq:LWbasic}), the process $(\Gamma\kgs{^{LW}}_{n\eta,nt})_t$ simplifies, noting also (\ref{alt_rep_mud}), to
    \begin{align}
		\Gamma^{\kgs{LW}}_{n\eta,nt}&= \left(\frac{2}{3}((n\eta)^3-n\eta)\sigma^2\right)^{-1/2}\nonumber\\
        &\quad\quad \times
		\left(\sum_{j=1}^{n\eta}w_\mu(\kgs{n\eta},j)X^{\kgs{LW}}_{ \lfloor nt\rfloor+j}-\sum_{j=1}^{n\eta}w_\mu(\kgs{n\eta},j)X^{\kgs{LW}}_{ \lfloor nt\rfloor- n\eta+j}\right), \nonumber
  \end{align}
\kgs{where $w_\mu(n\eta,j):=2j-n\eta-1$. If we combine the sums due to equality of the factors $w_\mu(n\eta,j)$ and replace $X^{LW}_{ \lfloor nt\rfloor+j} = (\lfloor nt\rfloor+j)\mu+b+\sigma Z_{\lfloor nt\rfloor+j}$ and analogously for $X^{LW}_{ \lfloor nt\rfloor- n\eta+j}$, this yields}
\begin{align}
        \Gamma^{\kgs{LW}}_{n\eta,nt}
        &\kgs{=\left(\frac{2}{3}((n\eta)^3-n\eta)\sigma^2\right)^{-1/2}}\nonumber\\
		&\kgs{ \quad \times
        n\eta\mu\sum_{j=1}^{n\eta}w_\mu(n\eta,j) + \sigma \sum_{j=1}^{n\eta}w_\mu(n\eta,j)\left(Z_{\lfloor nt\rfloor+j}-Z_{\lfloor nt\rfloor-n\eta+j}\right) }\nonumber\\
		&=\left(\frac{2}{3}((n\eta)^3-n\eta)\right)^{-1/2}  \sum_{j=1}^{n\eta}w_\mu(h,j)\left(Z_{\lfloor nt\rfloor+j}-Z_{\lfloor nt\rfloor-n\eta+j}\right),\label{rep_gamma_null}
	\end{align}
    where we used that $\sum_{j=1}^{n\eta}w_\mu(n\eta,j)=0$.
    
    As the scaling was chosen as the square root of the variance of the nominator,  
    for fixed $t$, $\Gamma^{\kgs{LW}}_{n\eta,nt}$ is standard normally distributed in $\mathbb R^2$. The covariance structure of each process component only depends on the window size $\eta$. \kgs{The limit behavior of $(\Gamma_{n\eta,nt}^{LW})_t$ is investigated in Section \ref{sect:konvergenz} with respect to a functional limit law.}\\

\subsubsection{\kgs{The statistical test and change point detection in practice}\label{sect:practice}}
    \kgs{Both for the LW and the RW, a} rejection threshold for the test can be obtained by simulation. In detail, for a given track $x_1,\ldots,x_T, T\in \mathbb N$ (cmp.~Figure \ref{fig:biv1} C), we derive the process $\left(G_{h,i}\right)_{i\in \{h,\ldots,T-h\}}$
    by shifting a double window of size $2h$ across the track to estimate the two dimensional expectation $\hat{\mu}(i-h,h)$ and $\hat\mu(i,h)$. The difference $\hat{\mu}\kgs{(i,h)}-\hat{\mu}\kgs{(i-h,h)}$ is scaled according to Definition \ref{def:ght}, resulting in the process $G$ whose two components are illustrated in Figure \ref{fig:biv1} D and E.  As a change in the expectation often affects both components, we can consider directly the two dimensional process $G$ (panel F) and use the maximal deviation $M_h(G):=\max_{i}\|G_{h,i}\|_{2}$ as a test statistic. Then we simulate independent realizations of the process $\left(\Gamma_{h,i}\right)_\kgs{i\in \{h,\ldots,T-h\}}$ and derive the rejection threshold $Q$ as the 95\%-quantile of the maximal deviations $M_h(\Gamma):=\max_i \|\Gamma_{h,i}\|_2$ of the simulated processes.  If $M_h(G)>Q$, i.e., if the process $\left(G_{h,i}\right)_\kgs{i\in \{h,\ldots,T-h\}}$ passes outside the circle with radius $Q$ around the origin, we reject the null hypothesis of no change points. For a rigorous justification neglecting rates of convergence see Remark \ref{remark}. 

    \kgs{Note that this testing procedure is consistent in the sense that the test statistic $M_h(G)$ tends to infinity under an alternative including at least one change point: We then have $M_h(G)=\max_{i}\|G_{h,i}\|_{2}\ge\|G_{h,i_0}\|_{2}$ where $i_0$ is the position of a change point. The nominator of (\ref{eq:ghtLW}) then has a nonzero limit as $h\to\infty$. Also the variance estimates in the denominator of (\ref{eq:ghtLW}) have finite limits. Hence, we obtain  $M_h(G)=\Omega(h^{3/2})$ tending to infinity as $h\to\infty$.} 
	
	After rejection of the null hypothesis, we estimate the number and location of change points using a procedure introduced in Messer et al. \cite{messerMFT,messerBivariate} (see Figure \ref{fig:biv1} H-K). Let $\hat{c}_{1}:=\text{argmax}_{i\in[h,T-h]}\|G_{h,i}\|_{2}$ be the estimator for the first change point. Then the process $G_{h,i}$ is deleted in the $2h$-neighbourhood around the change point $[\hat{c}_{1}-h+1,\hat{c}_{1}+h]$ because in that area, $G$ may be affected by the change point. Then we successively identify the maximal deviation of the remaining process $G$ to add change point estimates to the set of estimated change points $\hat{C}=\{\hat{c}_{1},\hat{c}_{2}\}$. This procedure is repeated until the remaining values of $G_{h,i}$ are inside the rejection area.

    \kgs{The computational complexity of the described algorithm is comparable for the RW and LW setting. For the first step, i.e., the statistical test procedure, the calculation of the test statistic $M$ for one track requires iteration over all time points in the range $[h,T-h]$, causing a factor of $(T-h)$ in the complexity, which we bound by $T$. As a second factor, the simulation of the rejection threshold $Q$ requires performing analogous computations on simulated tracks in $S$ simulations, causing a factor of $S$, which results in a complexity of $\mathrm{O}(T\cdot S)$ for the statistical test. In the second step, the change point detection algorithm iterates a procedure of complexity $\mathrm{O}(T)$ until no further change points are detected, resulting in a computational complexity of $\mathrm{O}(T\cdot k)$ for the change point detection, where $k$ denotes the number of change points.}

	\subsection{The Limit Process}\label{sect:konvergenz}	
	In this section we show the convergence in distribution  of the process of differences $(\Gamma^{\kgs{LW}}_{n\eta,nt}):=(\Gamma^{\kgs{LW},(1)}_{n\eta,nt},\Gamma^{\kgs{LW},(2)}_{n\eta,nt})$ from Definition \ref{def:ght} towards a Gaussian limit process whose distribution does not depend on the model parameters. Note that by independence of the components of $(\Gamma^{\kgs{LW},(1)}_{n\eta,nt}, \Gamma^{\kgs{LW},(2)}_{n\eta,nt})$, it is sufficient to show process convergence of the components individually.
 
	Recall that we fix $T>0$ and $0<\eta<\frac{T}{2}$. We denote $\para:=[\eta,T-\eta]$. We obtain the limit process with independent coordinates, which are distributed  as the centered Gaussian process $\Upsilon=(\Upsilon_t)_{t\in\para}$ with the covariance function 
	\begin{align}
		C:\para\times \para\to\Rset,\quad (s,t)&\mapsto \kappa\left(\frac{|s-t|}{\eta}\right),\label{eq:limit_covariance}
	\end{align}
	where the normalized autocovariance function $\kappa: [0,\infty)\to \Rset$ is given by (cmp.~also Figure \ref{fig:biv1} G)
	\begin{align}\label{eq:covariance_kappa}
		\kappa(x):= \left\{ \begin{array}{cl}
			3x^3-3x^2-\frac{3}{2}x+1, & \mbox{for } 0\le x \le 1, \vspace{2mm}\\
			-x^3+3x^2-\frac{3}{2}x-1, & \mbox {for } 1\le x \le 2,\vspace{2mm}\\
			0, &\mbox {for } x \ge 2.
		\end{array}\right.
	\end{align}
	
	We now consider the first component $d=1$ and require to normalize time. We set, cf.~\eqref{eq:GammaLW},
	\begin{align}\label{eq:upsiloncontaininggu1}
		\Upsilon^{(n)}_t &:=\Gamma^{\kgs{LW},(1)}_{n\eta,nt}, \quad t\in\para,
	\end{align}
	where $\Upsilon^{(n)}:=(\Upsilon^{(n)}_t)_{t\in\para}$ are considered as processes in the c\`{a}dl\`{a}g space $D(\para)$ endowed with the topology induced by the supremum norm $\|\,\cdot\,\|_\infty$. 
	
	The covariance functions $C_n:\para\times \para\to\Rset$ of the  processes $\Upsilon^{(n)}$ converge uniformly to the covariance function $C$ of the limit process, see (\ref{eq:limit_covariance}), as stated in the following Proposition.
	
    \begin{proposition}\label{lemma:convergencecovariancefunction} We have uniformly in $s,t\in\para$ that
		\begin{align*}
			C_n(s,t)\to C(s,t) \quad (n\to\infty).
		\end{align*}
    \end{proposition}
	
	The proof can be found in Appendix \ref{section:covariancefunctionupsilon}. The idea is to use that the process $\Upsilon^{(n)}$ equals in distribution a weighted sum of independent standard normal variables $\{\kgs{V}_j\,|\,j\in\Zset\}$, cf.~also (\ref{rep_gamma_null}),
	\begin{align}
		\Upsilon^{(n)}_t 
		\stackrel{d}{=}\sqrt{\frac{3}{2}}\left((n\eta)^{3}-n\eta\right)^{-1/2}\sum_{j=-\infty}^{\infty} g^{(n)}_t(j)\kgs{V}_j,\label{eq:upsiloncontaininggu2}
	\end{align}
	where, for $u \in\Rset$, $g^{(n)}_u$ is a function, see \eqref{def_gnu}, containing the weights in the estimator $\hat{\mu}$. Then the covariance function $C_n(s,t)$ can be interpreted as a Riemann sum, which uniformly converges to the function $C(s,t)$ as in equation \eqref{eq:limit_covariance}. The proof can be found in Appendix \ref{section:covariancefunctionupsilon}.
	
    Then as main result in this section we have the following process convergence:
	
    \begin{proposition}\label{proposition:convergenceupsilonn} For the processes $\Upsilon^{(n)}$ and $\Upsilon$ defined above we have convergence in distribution
		\begin{align*}
			\Upsilon^{(n)} \stackrel{d}{\longrightarrow} \Upsilon \quad (n\to\infty)
		\end{align*}
		within the c\`adl\`ag space $(D(\para),\|\,\cdot\,\|_\infty)$.
	\end{proposition}
	
	The proof of Proposition \ref{proposition:convergenceupsilonn} consists of showing convergence of the finite dimensional marginals (fdd convergence) and tightness. We use a criterion for the space  $(D(\para),\|\,\cdot\,\|_\infty)$, which can be found in Pollard \cite[Section V, Theorem 3]{Pollard1984}. The fdd convergence follows from the convergence of the covariance functions in Proposition \ref{lemma:convergencecovariancefunction} and L\'evy's continuity theorem. To verify the condition for tightness we exploit that the $\Upsilon^{(n)}$ are Gaussian processes, allowing to apply a bound of Dirksen \cite{dir15}, which is based on Talagrand's \cite{ta05,ta21} generic chaining method. The proof of Proposition \ref{proposition:convergenceupsilonn} can be found in Appendix \ref{app_func_proof}.\\
	
	\begin{remark}\label{remark}
	    Note that the convergence in Proposition \ref{proposition:convergenceupsilonn} can also be obtained for a time continuous version $G^{(n),\kgs{LW}}$ of $G_{h,i}^{LW}$ from  equation \eqref{eq:ghtLW} defined analogously  to the process $(\Gamma^{\kgs{LW}}_{n\eta, nt})$ in Definition \ref{def:ght}. This requires to strengthen the consistency in Proposition \ref{lemma:consistencyhatsigmasquared} to uniform consistency, so that the Lemma of Slutsky can be applied.
	\end{remark}

    \subsection{Simulations and multi scale approach}\label{sec:bivariatesimulations}
\begin{figure}[!h]
		\centering
        \includegraphics[width=1\textwidth]{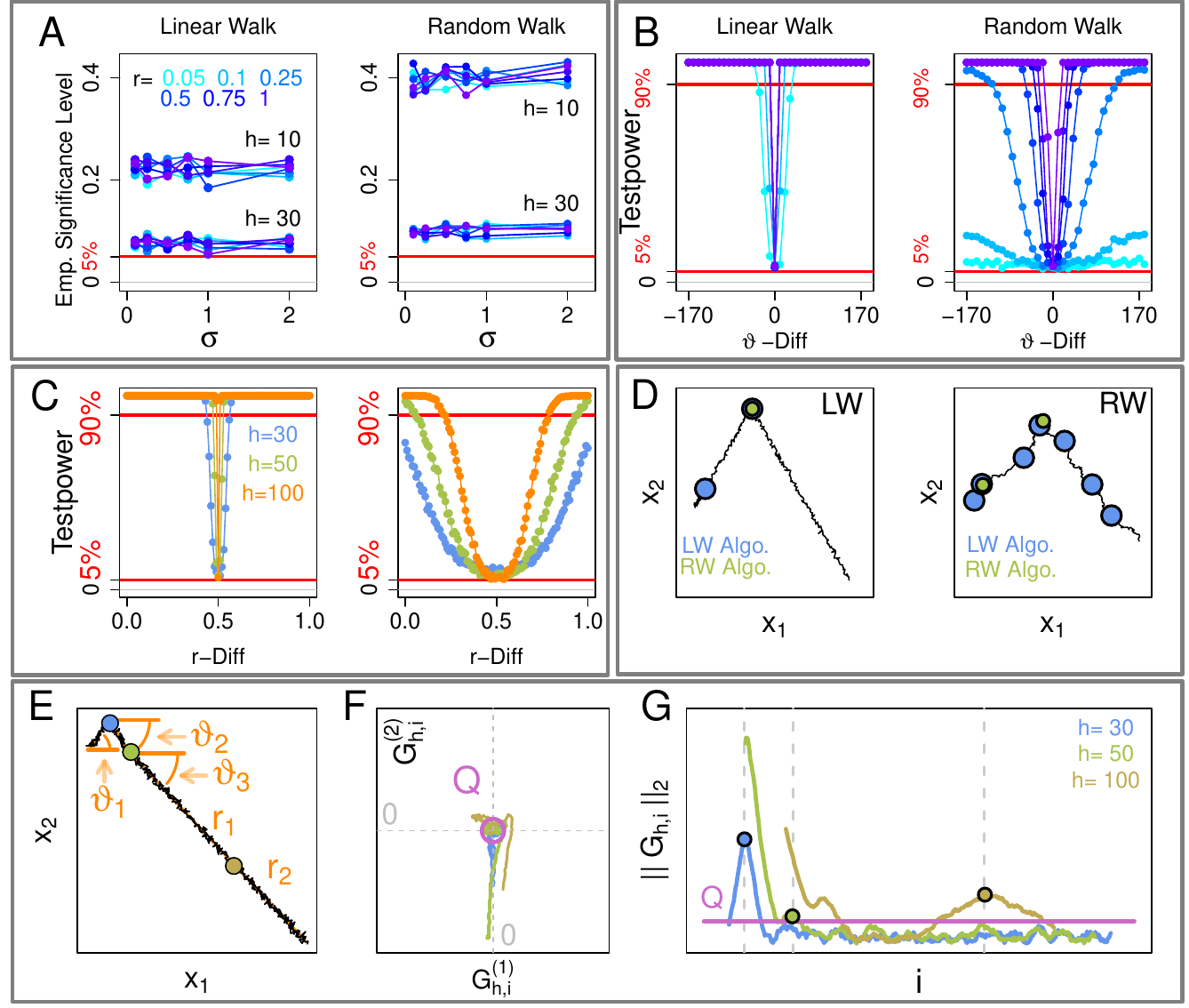}
		\caption{\small{(A--D) Simulations on the bivariate CPD comparing LW and RW, $T=400, \vartheta=\kgs{35/180\cdot\pi}, b=(0,0), \sigma^2=0.5,  1000$ simulations per parameter set. (A) Empirical significance level as a function of $\sigma^2$ and step length (color coded), $h\in \{10, 30\}$. (B) Test power on tracks with one change point in $\vartheta$ for varying magnitude of this change and varying step length, change point at $t=200$. (C) Test power on tracks with one change  in step length for varying magnitude of this change and varying window sizes (color coded). Other parameters are: $r_{pre-cp}=0.5, h\in\{30,50,100\}$. (D) Application of the bivariate CPD from both models to a LW (left) and a RW (right). Parameters for both tracks and CPD algorithms are $\vartheta=\kgs{(55,55,-55)/180\cdot\pi}$, $r=(0.2,1,1)$, $b=(0,0)$, $\sigma^2=1$, $T=300$, change point set $\{80, 150\}, h=30$. (E\kgs{-G}) Multi window CPD. (E) Simulated LW track with two changes in direction and one change in step length, true change points marked as points, $T=530$, $\vartheta=\kgs{(55, -55, -45, -45)/180\cdot\pi}$, $r=(1, 1, 1, 0.85)$, $\sigma=3$, $b=(0, 0)$, $CPs = (50, 110, 345)$. \kgs{(F) Process $G_{h,i}$ based on the track in E for window sizes in $\mathcal H=\{30, 50, 100\}$, colours as indicated in G. (G)} The absolute deviation of the difference processes $G_{h,\kgs{i}}$ for different window sizes $h$ as a function of time. Colours of window sizes indicated \kgs{in the top right corner}. Vertical lines indicate true change points, circles indicate estimated change points.}}
		\label{fig:biv5}
	\end{figure}

	\paragraph{\kgs{Significance level and comparison of LW and RW}}
	Because the results in the previous section are asymptotic, it is important to investigate the empirical significance level for finite data sets. Here we compare simulations of the classical MOSUM approach in the RW setting described in Section \ref{sec:MOSUMRW} with the new moving kernel approach in the LW setting described in Section \ref{sec:MOSUMLW}. Figure \ref{fig:biv5} A indicates that a significance level of 5\% is obtained for a window of about $h=30$ for the LW, while it is still considerably increased for $h=30$ in the RW  due to the summation of error terms in the test statistic and the higher variability in the underlying  tracks. We therefore recommend to choose at least $h=30$ for the LW  and $h=50$ for the RW. Similarly, due to the smaller variance, the test power is higher for the LW than for the RW case both for changes in direction and in step length (Figure \ref{fig:biv5} B and C).
	
	Due to the difference in model assumptions, erroneous application of the RW method to LW tracks typically is  conservative and yields an underestimation of change points (Figure \ref{fig:biv5} D, left). This is because the RW shows more variability in positions due to the summation of error terms, and this variability is inherited by the process $\left(G_{h,t}\right)_t^{(RW)}$. Therefore the rejection threshold in the RW case is generally higher. Vice versa,  erroneous application of the LW to RW tracks tends to yield a higher number of falsely detected change points (Figure \ref{fig:biv5} D, right).

\begin{figure}[htbp]
		\centering
		\includegraphics[width=1\textwidth]{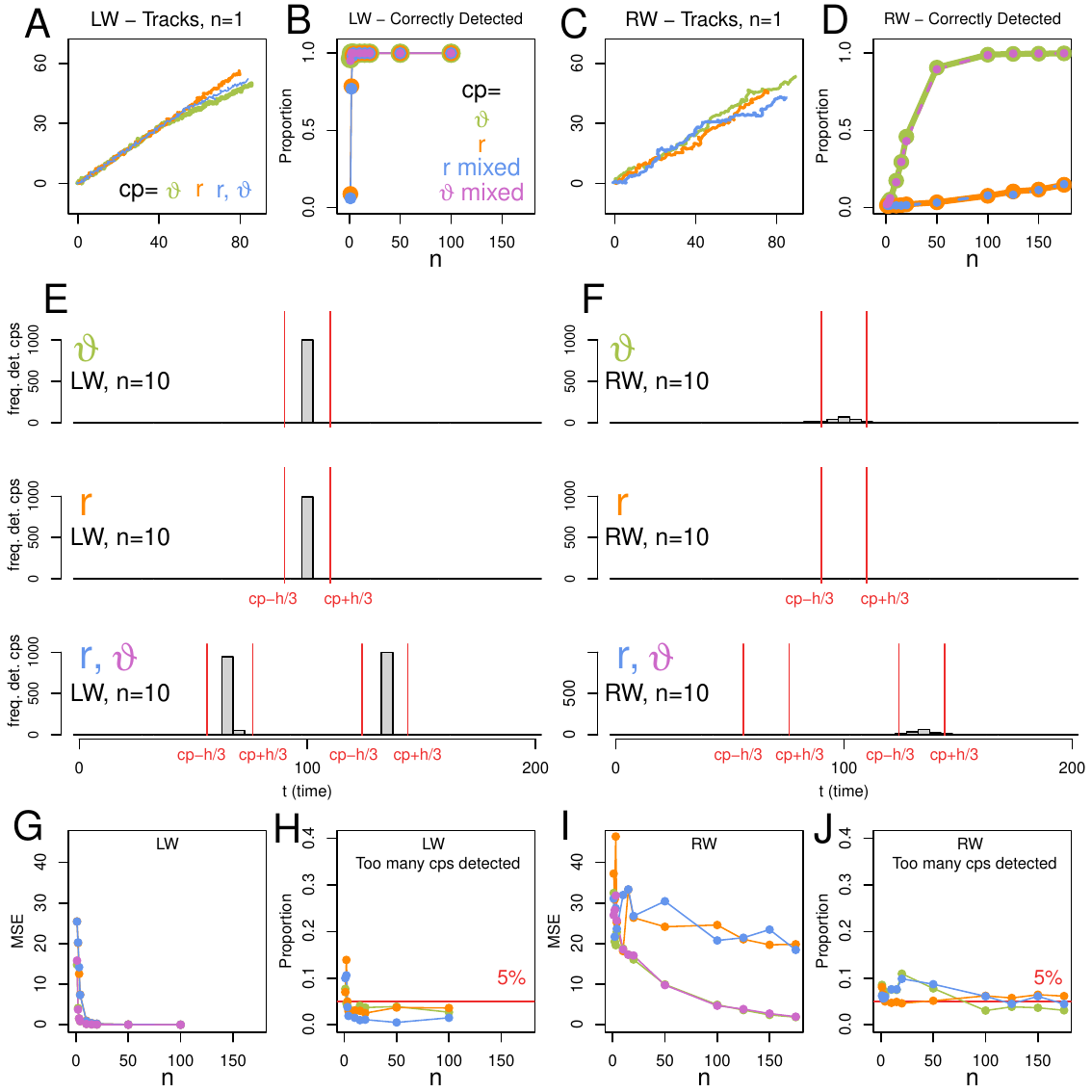}
		\caption{\small{\kgs{Asymptotic behavior of statistical test and change point detection for LW (A, B, E, G, H) and RW (C, D, F, I, J) in three scenarios, simulated with parameters $b=(0,0)$, $\sigma=0.5$, $T=200$. The three scenarios are (example tracks in A, C): (1, green) one change point at time $t = 100$ in direction with $\vartheta=(35/180)\pi$ before and $\vartheta=(25/180)\pi$ after the change point, (2, orange) one change point in step length at time $t = 100$  with $r$ changing from $0.5$ to $0.48$, and (3, blue) a combination of the same two change points, $r$ changing from $0.48$ to $0.5$ at $t=66$ and $\vartheta$ changing from $(25/180)\pi$ to $(35/180)\pi$ at $t=134$. Performance of the new kernel approach was analysed in $1000$ simulations in the LW and of the classical MOSUM in $2000$ simulations in the RW using $\eta=30$ and different values of $n$. For the asymptotic setting we let time $\lfloor nt\rfloor$ and window size $n\eta$ grow linearly in $n$. Proportion of simulations with correctly detected change points in the LW (B) and RW (D). Histograms of estimated change points for all three scenarios for $n=10$ in the LW (E) and RW (F), where the time points of detected change points are scaled with $1/n$ for better comparability.  A true change point was then called correctly detected if an estimated change point fell within a distance of $h/3$ of the true change point. Mean squared error (MSE) of correctly detected change points in the LW (G) and RW (I). Proportion of simulations in which too many change points were detected in the LW (H) and the RW (J).}}}
		\label{fig:biv6}
	\end{figure}

\kgs{\paragraph{Number of falsely detected change points under the null hypothesis and in the presence of true change points} Note that although the algorithm for change point detection seems to implicitly apply multiple testing for the successive identification of change points, this does not affect the probability of falsely identifying change points both in the LW and RW \citep[see also ][]{messerMFT}. First, under the null hypothesis of no change points, the statistical test keeps the significance level (Figure \ref{fig:biv5} A), and in the rare case of falsely detected change points, only few falsely detected change points were observed in our simulations. The proportion of simulations with two falsely detected change points under the null hypothesis was $1.6$ \% across all simulations in the LW case ($0.15$ \% for RW), while more than two change points were falsely detected only in $0.9$\% of simulations  ($0.05$ \% for RW). Second, if a change point exists, it only affects the process $G$ in its $2h$-neighborhood, which is cut out after detection of the change point. Therefore, existing change points should not affect the behavior of the procedure outside their $2h$-neighborhood and thus, not enhance the probability of falsely detecting additional change points. In Figure \ref{fig:biv6} we illustrate the asymptotic performance of the proposed procedures, i.e., as $n$ increases. Panels H and J show the proportion of simulations in which too many change points were detected, indicating that as $n$ increases, the probability of falsely detecting additional change points is not increased also in the presence of true change points.}

\paragraph{\kgs{Precision of change point detection}}
\kgs{Figure \ref{fig:biv6} also illustrates simulations of three different change point scenarios in $\vartheta$ and $r$ in which the mean squared error of the estimated change point tends to zero and the probability of correctly detecting an existent change point tends to one as $n$ increases. The latter is in agreement with the divergence of the test statistic at a change point discussed in Section \ref{sect:practice} and is illustrated both for the LW and the RW case. Note that consistency of the number and location of change points for the RW in the classical MOSUM setting has been shown by \citep{kirch2021moving}. However, the results for the RW case do not directly extend to our LW case due to the prefactors appearing in our kernel estimator. We leave this issue for further research. }

\paragraph{\kgs{Multi scale approach}}
    Interestingly, the present methods are compatible with the combined use of multiple windows $\kgs{\mathcal{H}}:=\{h_{1},h_{2},\ldots,h_{\kgs{H}}\}, \kgs{H}\in\mathbb{N}$   \citep[cmp.~][]{messerMFT,messerBivariate} \kgs{both for the RW and the LW}. To that end, one can calculate the process $G_{h_\kgs{j},i}$ for every window size  $h_\kgs{j}\in \kgs{\mathcal{H}}$ and then use the maximum $M:=\max_{\kgs{j}=1,\ldots,\kgs{H}}M_{h_\kgs{j}}$ of all maxima across processes as a test statistic\kgs{, where $M_{h_\kgs{j}}:= \max_{i} \|G_{h_\kgs{j},i}\|^2$}. The threshold $Q$ can be simulated analogously to the one window approach. That means we simulate, for every $h_\kgs{j}$, the process $\Gamma_{h_\kgs{j},i}$ (Definition \ref{def:ght}), where for each individual simulation, the processes of all window sizes $h_\kgs{j} \in \kgs{\mathcal{H}}$ are based on the same realization of the random sequence $\textbf{Z}:=\left(Z_i\right)_{i\in\mathbb{N}_0}$. Thus, the processes $\Gamma_{h_\kgs{j},i}$ are related across the different windows in the same way as the processes $G_{h_\kgs{j},i}$. We then obtain the distribution of the maximum $M$ of all maximal distances of $\Gamma_{h_\kgs{j},i}$ to the origin. Its $95\%$-quantile serves as the new rejection threshold $Q$.
	
	The CPD is performed by using the global threshold $Q$ on each process $\left(G_{h_\kgs{j},i}\right)$. The sets of change points derived for each individual window are then combined into a final set starting with the change points estimated in the smallest window and then successively adding change points from larger windows to a set of 'accepted' change points if their respective $2h_\kgs{j}$-neighborhood does not contain an already accepted change point \citep[cmp.~][]{messerMFT}.

    \kgs{Application of multiple windows then increases the complexity by a factor of $H$, where $H$ denotes the number of applied windows. That means, the computational complexity of the change point test is $\mathrm{O}(H\cdot T\cdot S)$, since every step needs to be performed for each of the $H$ window sizes. The same holds for the change point detection algorithm, resulting in a computational complexity of $\mathrm{O}(H\cdot T\cdot k)$.}

	Figure \ref{fig:biv5} E illustrates an example in which two large change points in direction in short succession are followed by one small change in step length. Panel F illustrates the processes $\kgs{G}$ for different window sizes, and panel G shows only the distances from the origin as a function of time. As one can see in panel G, only the combination of multiple window sizes allowed the detection of all change points. While short term changes could be detected by smaller windows, small change on a large time scale required a larger temporal window. Thus, a combination of multiple windows is particularly interesting if one aims at detecting changes at multiple time scales. \kgs{As the computational complexity increases with the number of windows, this number should however remain small. While the smallest window should be sufficiently large to maintain the asymptotic significance level, larger windows should be sufficiently small to distinguish between multiple change points. In addition, because all processes $\left(G_{h_\kgs{j},i}\right)$ are based on the same underlying track, they show high dependency for similar window sizes, and test power therefore hardly increases when adding highly similar window sizes to the window set. In order to reduce computational complexity, chosen window sizes should therefore be sufficiently different from each other \citep[cf.~][]{messerMFT}.}

    \section{Application} \label{sec:application}

    Here we apply the statistical test and CPD to example tracks from the data set of organelle movements introduced in Section \ref{sec:motivation}. The data set consisted of $41$ planar tracks representing the movement of cell-organelles of the type plastid and of $12$ planar tracks of peroxisomes in root cells of \textit{Arabidopsis thaliana}. 
   
    We first present the analysis in the LW model and show a comparison to the detected change points within the RW in section \ref{sect:LWRW}. We use a window size of $h=30$ to adhere to the chosen significance level of $5\%$ as suggested by the simulations in Section \ref{sec:bivariatesimulations}. As we will see, a high number of big changes is estimated with this single window only, and estimated change points are sometimes quite close together. We therefore used only the smallest possible window instead of a set of multiple windows in order to keep the test power of this individual window as high as possible.
	
\begin{figure}[htbp]	
		\centering
        \includegraphics[width=1\textwidth]{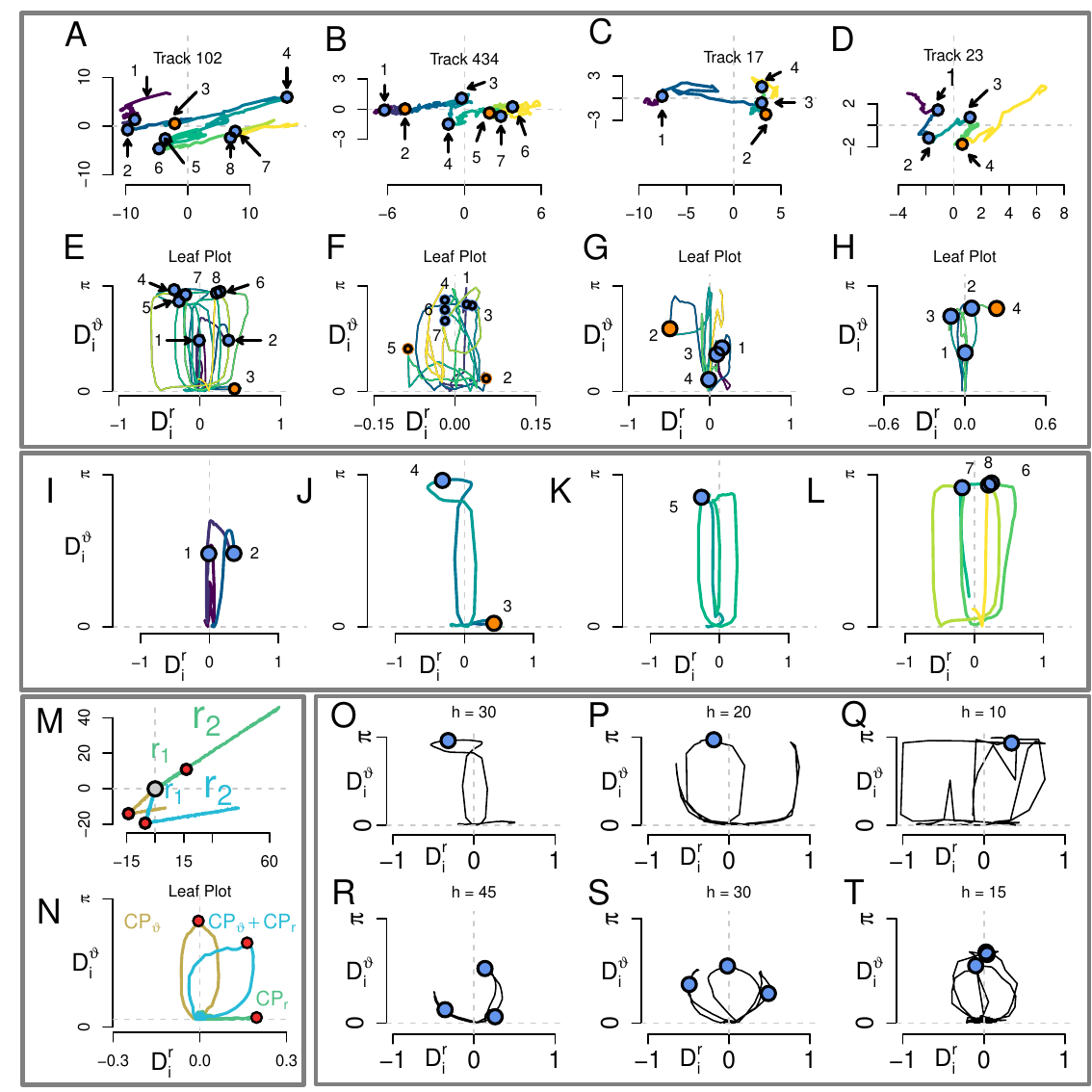}
		\caption{(A, B) Tracks 102 and 434 from the plastid data set. Change points estimated with the $CPD_{LW}$ and $h=30$ (blue+orange dots). Tracks between change points color coded. (C)+(D) Tracks 17 and 23 from the peroxisome data set analysed with $CPD_{LW}$ and $h=30$. (E)-(H) Corresponding leaf plots to (A)-(D), created with $h=30$ in (E, F, H) and $h=20$ in (G). (I-L) Leaf plot from (A) separated in different change points. (M) Simulated LW tracks with Parameters (pale brown track) $\vartheta=\kgs{(-135,10)/180\cdot\pi}$, $r=0.1$, $b=(0,0)$, $\sigma=0.15$, $T=400$, $t_{CP}=200$, (green track) $\vartheta=\kgs{(35,35)/180\cdot\pi}$, $r=(0.1,0.5)$, $b=(0,0)$, $\sigma=0.25$, $T=300$, $t_{CP}=100$, (blue track) $\vartheta=\kgs{(-105,10)/180\cdot\pi}$, $r=(0.1,0.25)$, $b=(0,0)$, $\sigma=0.25$, $T=400$, $t_{CP}=200$ with true change points (red). (N) Leaf plot of the tracks with the corresponding colour from (M). (O-Q) Leaf plots of the loop of change point no 4 in Track 102 with $h=30$ (O), $h=20$ (P) and $h=10$ (Q). (R-T): Leaf plots of a simulated track with small standard deviation and three change points analyzed with $h=15$ (T), $h=30$ (S) and $h=45$ (R). Simulated leaves are separated for small window sizes in T and S and merge for bigger windows (R).}
		\label{fig:appli9}
	\end{figure}

    Figure \ref{fig:appli9} A and B show our two example tracks, where blue and orange points indicate change points in direction and step size, respectively, estimated within the LW model. As one can see, the LW method estimates  a number of change points, particularly at prominent changes in the movement direction, that are visible by eye.

    \subsubsection*{The leaf plot}\label{sect:leafplot} 
	In order to analyse the detected change points in more detail and to distinguish between change points in the direction and in the step length, we propose to use the original parameterization of direction and step length in a graphical approach. To that end, we calculate the process of differences in direction, $D_\kgs{i}^{\vartheta}:=D_\kgs{i}(\hat{\vartheta}_{\ell}\kgs{(i)},\hat{\vartheta}_{r}\kgs{(i)})$ and in step length, $D_\kgs{i}^{r}:=D_\kgs{i}(\hat{r}_{\ell}\kgs{(i)},\hat{r}_{r}\kgs{(i)})$ in a moving double window as described in equation \eqref{eq:mosum}. We then plot these difference processes against each other in a bivariate plot.
	
	This plot will be called \textit{leaf plot} for short because ideal change points show up in leaf form in the bivariate difference process representation. Figure \ref{fig:appli9} M and N show simple examples of such leaves for three simulated LWs  with different change points. The green process shows a change point in  step length, which is illustrated in the leaf plot as a degenerated flat leaf in which only the difference process $D_t^{r}$ on the $x$-axis deviates from zero. Second, in pale brown, we see a track with a change only in the movement direction in panel M. The corresponding leaf in panel N deviates from the origin in the vertical direction, where the maximal deviation is expected at the time of the change point in direction. Note that this leaf also shows slight horizontal deviations around the change point in direction. This does not indicate changes in  the step length but is caused by the fact that a change point in the direction  tends to slightly bias the estimates in the step length in the area around this change point (cmp.~Supplement 1). Third, in light blue, we see a change  in direction and step length occurring at the same time, which shows up in a diagonally oriented leaf in panel N. \\
	
	For track 102, we observe the leaf plot in panel E, where its leaves, i.e., the sections of the process around the estimated change points, are depicted individually in panels I-L. The numbers next to the blue and orange points indicate the temporal order of detected change points in the track. Interestingly, most estimated change points show up in vertical deviations in the leaf plot, with a maximal difference in directions of almost  $\pi$, suggesting strong changes in the movement direction (blue points). 
	
	In addition to these strong changes in the movement direction, we also observe one horizontal deviation in the leaf plot at change point no.~3 (orange point), suggesting an increase in step length at this location.
	
	In order to interpret leaves that do not correspond completely to the stereotypical form indicated in panel N, we observe that such phenomena will occur if successive change points are closer than the window size $h$. This is illustrated in Figure \ref{fig:appli9} O-T by comparison to a simulation. The respective loop of track 102 is shown in panel O with an analysis window $h=30$. Decreasing the window size to $h=20$ (panel P) yields one vertical leave indicating a change point in direction and another structure indicating additional change points on even smaller time scales. Reducing the window size further to $h=10$ (panel Q) however does not reveal new leaves due to high empirical variability. Therefore, panels R-T illustrate a simulation with small variability and comparable effects. In panel T, we observe three distinct leaves corresponding to three simulated change points (blue points) in short succession in a leaf plot with small $h=15$. Successively increasing the window size yields similar structures as observed in the track (panels S and R).
	
	Therefore, the leaf plot may allow the visual classification of detected change points into directional changes (vertical leaves), step length changes (horizontal leaves) and simultaneous changes (diagonal leaves). In addition, loop like structures can be indicative of the analysis window overlapping more than one change point.

	\subsubsection*{Further examples}
	As another example of a plastid track, Figure \ref{fig:appli9} B and F shows the estimated change points of Track 434. The results are highly similar, showing a number of strong changes in direction by about $\pi$ (e.g., change points no.~3, 4, 6, 7, blue points) as well as changes in the step length (change points no.~2 and 5, orange points). 		Panels \ref{fig:appli9} C, D, G, H show further examples of peroxisome tracks, which indicate similar behavior. Most interestingly, the peroxisome tracks show clear changes in the step length. For example, track 17 shows a higher step length between change points no.~1 and 2 (Figure \ref{fig:appli9} C, orange point) indicated by the horizontal leaves around change point 2 (panel G, orange point). Similarly, peroxisome track 23 (panels D and H) shows a strong change in the speed at the last change point no.~4 (orange point in panel E), as indicated by the horizontal deviation on the right of the leaf plot (panel H, $h=30$).

	\subsection{Comparison of the LW and RW approach}\label{sect:LWRW} 
\begin{figure}[htbp]
		\centering
        \includegraphics[width=1\textwidth]{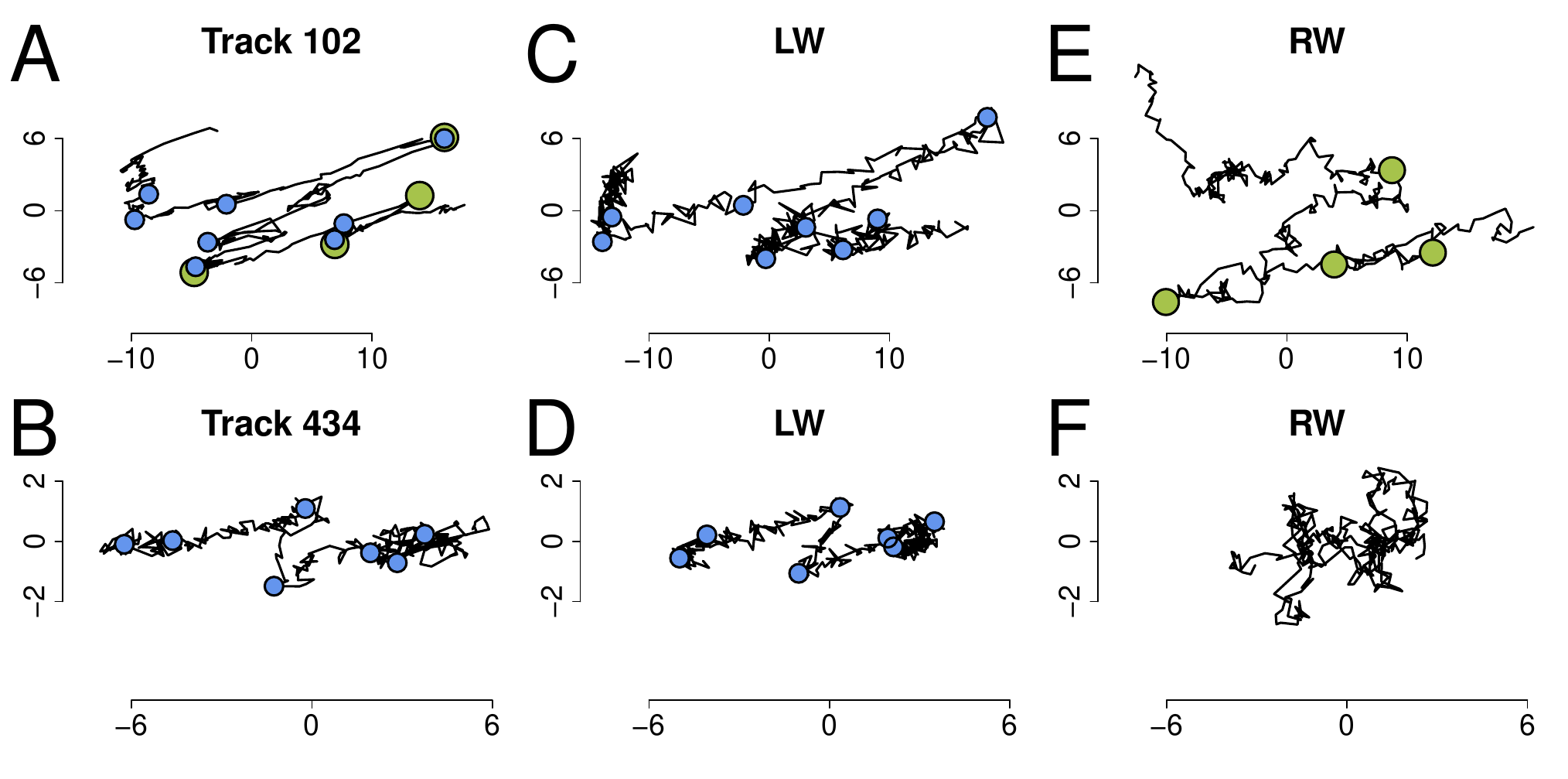}
		\caption{Tracks 102 (A) and 434 (B) with change points estimated by with the $CPD_{LW}$ (blue dots) and $CPD_{RW}$ (green dots) marked, h=30. Simulated LW-tracks with the estimated $CPD_{LW}$ change points, direction and step length estimated between change points, variance estimated as median \kgs{across time points $i\in\{1,\ldots,T-2h\}$ of} $\kgs{(\hat\sigma^{LW})^2(i,2h)}$\kgs{, where the values in a $2h$-neighbourhood around each detected CP were removed, based on} track 102 (C) and track 434 (D). Simulated RW-tracks with the estimated $CPD_{RW}$ change points, direction and step length estimated between change points, variance estimated as median \kgs{across time points $i\in\{1,\ldots,T-2h\}$ of} $\kgs{(\hat\sigma^{RW})^2(i,2h)}$\kgs{, where the values in a $2h$-neighbourhood around each detected CP were removed, based on} track 102 (E) and track 434 (F).}
		\label{fig:appli3}
	\end{figure}

	As explained above in Section \ref{sect:Null} and \ref{sect:discmodelass}, the LW was explicitly designed to describe movements along roughly linear cellular structures. Accordingly, the change points detected with the LW method often agree with visual inspection, suggesting a number of strong changes in movement direction, and also changes in step length. In comparison, the RW method estimates considerably fewer change points. In track 434, the null hypothesis of no change points was not even rejected. This is illustrated in Figure \ref{fig:appli3} B, which shows no green circles, i.e., no estimated change points in the RW approach. For track 102, the RW approach estimated only four change points (green circles in Figure \ref{fig:appli3} A), particularly missing visually prominent changes in direction. This is in line with the observation that the RW-CPD tends to overlook change points when applied to an LW (Section \ref{sec:bivariatesimulations}).

	As an additional visual comparison, we simulated LW and RW tracks with piecewise constant movement direction and speed in Figure \ref{fig:appli3} C-F. To that end, we estimated the piecewise constant movement direction and step length from the organelle tracks within each section between  detected change points, both for the LW and the RW, using the model specific estimators of parameters and change points described earlier. We then simulated LW and RW tracks with these estimates. Again, the LW shows strong linear components similar to the shown organelle tracks, while the RW has weaker linear parts and a tendency to show a curled behavior.
	
\section{Discussion}
	This paper followed three main aims. The first goal was to derive a stochastic model with which movements of cell organelles can be described. Specifically, a data set of plastid movements indicated that organelles can move in a rather linear manner, following piecewise linear structures with abrupt changes in movement direction and speed. Therefore, the second aim of the paper was to develop an analysis method that can statistically test for changes in the movement direction and speed. Third, if the null hypothesis of no changes is rejected, we aimed at estimating the locations of the change points and at distinguishing between changes in the direction and the step length.
	
	Concerning the first aim, we have presented a new stochastic model for the description of two dimensional movement patterns. The model is called here a linear walk (LW) in order to stress similarities as well as differences to the \kgs{biased} random walk (RW) assumptions. On the one hand, the LW assumes, similar to the RW, that movement has a piecewise constant direction and speed and that each step is drawn from a probability distribution. On the other hand, movement within a LW follows more strictly a straight line than the RW by assuming that observations are equidistant points on a straight line offset by independent random errors. As a consequence,  increments in a LW show a clear dependence structure of order one and allow for a relatively strict movement along linear structures as is observed in  cell organelle movements reported here.
	
	In order to address the second aim, i.e., to test for changes in movement direction, we used a moving window approach and first  explained parameter estimation within the LW. Illustrating that a univariate analysis of the movement direction requires assuming a constant step length, we switched to a bivariate setting  in which the polar parametrization is transformed to a cartesian parametrization. \kgs{A classical MOSUM approach could be applied in the RW case but is }limited in the LW setting. We therefore used the new MLEs to derive a  new moving kernel statistic replacing the MOSUM process. We then showed convergence of the resulting process, assuming the true scaling, to a limit Gaussian process which is independent from  the model parameters, and we also showed pointwise strong consistency of the variance estimator. These theoretical results motivated the use of our new moving kernel statistic within the LW, where we replaced the true scaling by the variance estimator in order to show good performance with respect to significance level and test power via simulations.
	
	Estimation of change points as described in the third aim can then be performed with a straightforward algorithm \kgs{both for the RW and LW}. In particular, we showed that the described method can be easily applied using a whole set of multiple window sizes in order to perform change point detection on multiple time scales simultaneously. \kgs{The complexity of the algorithm is linear in the length of the time series, the number of simulations used for the derivation of the rejection threshold, in the number of windows and the number of change points.} Finally, we proposed a graphical technique called the leaf plot in order to distinguish qualitatively between change points in the direction and step length. In the leaf plot, local differences in the estimated step length are plotted against local differences in the estimated movement direction in a moving window analysis, thus reverting back to the polar parameterization. The form of the resulting leaves allows a qualitative discrimination between changes in the individual parameters. The analysis of the sample data set indicated a high number of estimated change points both in the direction and step length, where most directional change points ranged around the value of $\pi$, indicating a frequent quasi reversal of the organelle movement direction.

\kgs{Concerning the model assumptions, both the LW and the (biased) RW used here share a number of advantages. Both models describe }visually observable sections of roughly constant movement direction. All model parameters have a clear and simple interpretation and can be estimated from the data. Also, by assuming absolute angles instead of the widely used turning angles, the LW \kgs{and  biased RW allow} to distinguish between phases of fast linear movement that show different absolute movement directions as observed in the present data set. In this respect, we consider such models more suitable here than, e.g., correlated random walks \citep[see also][]{Plomer2024b}.
 
\kgs{In the comparison between the LW and RW, we observe (see Fig.~\ref{fig:Discussionmodelassumption}) that the LW sticks quite closely to a straight line, while the RW shows a much higher variability. As a consequence, change point methods based on the RW have a lower test power and could fail to detect visually prominent changes in the movement direction in the present data set, while the change points estimated with the LW tended to correspond more closely to visual inspection.}
	
	However, one should also note that the visual correspondence of the LW with the patterns of organelle tracks has limits. In particular, LW tracks do not appear perfectly organic due to the negative serial correlation of increments, which is caused by the independence of offsets of each measurement point from the line. 
    It would be interesting to study extended models in which the strong serial correlation is weakened, for example by a slow backward drift towards the underlying line, \kgs{or by a combination of RW assumptions with a stricter movement direction}. However, this often requires further technical assumptions such as for example attraction points, whose number and location typically needs to be pre-specified.
    Such assumptions would, therefore, come at the cost of higher technical complexity. 

	Considering the statistical test and change point detection algorithm, we have shown that methods for multiple window sizes derived earlier \citep{messerMFT} can easily be applied in the present setting. Thus, the present analysis also allows for the analysis and estimation of change points in movement direction and speed on multiple time scales. One should, however, note that when applying the asymptotic procedure, a minimal window of size $h=30$ (for the LW) should be applied in order to keep the required significance level. As a consequence, change points that occur in faster succession cannot be detected. They will, however, show up in loops and other shapes in the leaf plot indicating the potential existence of additional change points.
		
	Finally, one should also note that although the statistical test and change point estimation were derived and presented for two dimensional movements, they are easily applicable to an arbitrary number of dimensions. This is because the theoretical results rely on the behavior of the processes in one single dimension, where errors in all directions are assumed independent. Therefore, these methods can be extended to an arbitrary number of dimensions. 
 
    In  Supplement 4, we provide an R code that implements the statistical test for and change point estimation as well as basic representation of a two-dimensional track with its estimated change points. The algorithm directly uses multiple windows for change point detection and illustrates the estimated change points in a leaf plot.

We thus believe that the present paper can contribute to the theoretical understanding and  application of change point analysis of biological movement patterns that show abrupt changes in their movement direction and speed.

\begin{appendix}

	\section{Proof of Proposition \ref{theorem:ml_mu_positions}}\label{proof:ml_mu_sigma_positions_t}
  	\kgs{Note that throughout the appendix all processes and estimators refer to the LW case. For ease of notation, we will therefore omit the superscript $LW$ in the following. Also, we will use the short notation $\hat\mu:=\hat\mu(i,h)$ and $\hat b:=\hat b(i,h)$.}
   
    Let $X_{i+1}, X_{i+2}, \ldots, X_{i+h}$ be a sequence of positions from a uni-directional\\ \textbf{LW}({$\vartheta,r,b$},$\sigma$,$\emptyset$). Then all $X_{i+j}^{(d)}$ are independent for all $j=1,\ldots,h$ and $d=1,2$ and with distribution
		\begin{align*}
			X_{i+j}^{(d)} \sim \mathcal{N}\left( (i+j)\mu^{(d)}+b^{(d)} , \sigma^2 \right),
   		\end{align*}
     with 
     \begin{align*}
     (\mu^{(1)},\mu^{(2)})=(r\cdot\cos(\vartheta),r\cdot\sin(\vartheta)).
		\end{align*}
		The corresponding log-likelihood is given by
		\begin{align}
			-h\log{2\pi\sigma^2}-\frac{1}{2\sigma^2}\sum_{d=1,2}\sum_{j=1}^{h} \left( X_{i+j}^{(d)}-(i+j)\mu^{(d)}-b^{(d)}\right)^2.\label{eq:likelihood_positions_t}
		\end{align}
		Differentiating with respect to $\mu^{(d)}$ and $b^{(d)}$ for fixed $d$ results in the following system of equations
		\begin{align}
			\sum_{j=1}^{h}(i+j)\left(X_{i+j}^{(d)}-(i+j)\hat\mu^{(d)}-\hat b^{(d)}\right)&=0, \label{eq:system1}\\
     \sum_{j=1}^{h}\left(X_{i+j}^{(d)}-(i+j)\hat\mu^{(d)}-\hat b^{(d)}\right)&=0.\label{eq:system2}
		\end{align}
        In \eqref{eq:system1} we use the short notation $s_1:=\sum_{j=1}^{h}(i+j)$ and $s_2:=\sum_{j=1}^{h}(i+j)^2$ to solve for $\hat\mu^{(d)}$ as a function of $\hat b^{(d)}$,
		\begin{align}
			&\hat \mu^{(d)}=\frac{1}{s_2}\sum_{j=1}^{h}(i+j)X_{i+j}^{(d)}-\frac{s_1}{s_2}\hat b^{(d)}.\label{eq:system3}
		\end{align}
        Now we solve \eqref{eq:system2} for $\hat b^{(d)}$ and replace $\hat \mu^{(d)}$ by the term in \eqref{eq:system3},
	       \begin{align*}
			\hat b^{(d)} 
            &= \frac{1}{h}\sum_{j=1}^{h}X_{i+j}^{(d)}-\frac{\hat\mu^{(d)}}{h}s_1\\
			&=\frac{1}{h}\sum_{j=1}^{h}\left(1-\frac{s_1}{s_2}(i+j) \right)X_{i+j}^{(d)}+\frac{s_1^2}{h\cdot s_2}\hat b^{(d)}.
		\end{align*}
        This yields
		\begin{align}
			& \hat b^{(d)}=\frac{1}{s_1^2-s_2 h}\sum_{j=1}^{h}\left(s_1(i+j)-s_2\right)X_{i+j}^{(d)}.\label{eq:system4}		
   \end{align}
  Inserting $s_1=h(h+2\kgs{i}+1)/2$ and $s_2=h(2h^2+g\kgs{i}h+3h+6\kgs{i}^2+6\kgs{i}+1)/6$ yields 
  \begin{align*}
  \hat b^{(d)}&=(h^3-h)^{-1}\sum_{j=1}^{h}w_{b}(h,i,j)X_{i+j}^{(d)}, \quad \text{where}\\
      w_{b}(h,i,j)&:=(-6hj-12ji-6j+4h^2+6hi+6h+6i+2).
  \end{align*} To derive $\hat\mu^{(d)}$ we plug in \eqref{eq:system4} in  \eqref{eq:system3}, which yields
		\begin{align}\label{alt_rep_mud}
			\hat\mu^{(d)}&=
			\frac{1}{s_2h-s_1^2}\sum_{j=1}^{h}\left(h(i+j)-s_1\right)X_{i+j}^{(d)}\nonumber\\
            &=\frac{6}{h^3-h}\sum_{j=1}^{h}(2j-h-1)X_{i+j}^{(d)}\\ 
            &= \frac{6}{h^3-h}\sum_{j=1}^{\lfloor\frac{h}{2}\rfloor}\left(h-(2j-1)\right)\left(X_{i+j}^{(d)}-X_{i+h-j+1}^{(d)}\right).\nonumber
		\end{align}
        For the estimation of $\sigma^2$ we differentiate \eqref{eq:likelihood_positions_t} with respect to $\sigma^2$ and obtain equation \eqref{eq:MLEsigmaLM}. For an interpretable notation of $\hat b$, we  note that 
        \begin{align*}
            \frac{w_{b}(h,i,j)}{6(2j-h-1)}= -(i+(h+1)/2)+ \frac{(h^2-1)}{6(2j-h-1)}.
        \end{align*}
        This implies
        \begin{align*}
            \hat b &= \frac{1}{h(h^2-1)}\sum_{j=1}^h{(h^2-1)X_{i+j}} - (i+(h+1)/2)\hat\mu\\ 
            &= \bar X_\kgs{i}- (i+(h+1)/2)\hat\mu,
        \end{align*}
        where $\bar{X}_\kgs{i}:=h^{-1}\sum_{j=1}^h{X_{i+j}}$ denotes the mean of the observed points. Unbiasedness of the MLEs is shown in Supplement 2.
        		
		\section{Proof of Part 1 of Proposition \ref{lemma:consistencyhatsigmasquared}}\label{appendix:proofstrongconsistencymub}\label{appendix:SLLNnew}
        To show the strong consistency of $\hat{\mu}^{(d)}$ we use a variation of the SLLN for a setting of weighted sums. The proof of the strong consistency of $\hat{\sigma}^2$ follows the same structure but involves more tedious calculations. Together with the second required SLLN and the proof for the strong consistency of $\hat b$, it can be  found in Supplement 3.
		\begin{proposition}[SLLN with weights]\label{proposition:SLLNnew}\quad\\				Let $X_j$ be independent with $\mathbb{E}[X_j]=0$, but not necessarily identically distributed random variables with $\sup_{j\ge 1}\mathbb{E}[X_j^4]< \infty$. Let $S_n:=\sum_{j=1}^{m_n} v(n,j)X_j$, where $m_n\in\mathbb{N}$ and $v(n,j)\in\mathbb{R}$ are weights. Then $\lim\limits_{n\rightarrow\infty}S_n=0$ almost surely if
				\begin{align}
					\sum_{j=1}^{m_n}v(n,j)^2 \kgs{=}  \mathrm{O}(n^{-1}). \label{eq:conditionSLLN1}
				\end{align}
		\end{proposition}
Proposition \ref{proposition:SLLNnew} can be proven by a standard argument based on the Lemma of Borel--Cantelli and Markov's inequality. 

				Strong consistency of $\mu^{(d)}$ is now shown by verifying  condition 1 in Proposition \ref{proposition:SLLNnew}. The MLE for $\mu^{(d)}$ in the time continuous setting is given by
			\begin{align*}
				\hat{\mu}^{(d)}(nt,n\eta)&=\frac{6}{(n\eta)^3-n\eta}\sum_{j=1}^{n\eta}\left(2j-n\eta-1\right)X^{(d)}_{\lfloor nt\rfloor+j}\\
                &=\sum_{j=1}^{n\eta}v(n,j)X^{(d)}_{\lfloor nt\rfloor+j}, \quad\text{with}\\
                v(n,j)&:=\frac{6\left(2j-n\eta-1\right)}{(n\eta)^3-n\eta}.
			\end{align*}
			Since $X^{(d)}_{\lfloor nt\rfloor+j}\sim \mathcal{N}\left((\lfloor nt\rfloor+j)\mu^{(d)},\sigma^2\right)$, the fourth moment of $X^{(d)}_{\lfloor nt\rfloor+j}$ exists and since $\hat{\mu}^{(d)}$ is unbiased, it suffices to show that as $n\rightarrow\infty$, $(\hat{\mu}(nt,n\eta)-\mu) \rightarrow 0$ almost surely. According to Proposition \ref{proposition:SLLNnew}, it is therefore sufficient to show that $\sum_{j=1}^{m_n}v(n,j)^2 \in \mathrm{O}(n^{-1})$. This holds because
			\begin{align*}
				\sum_{j=1}^{n\eta}v(n,j)^2 &= \frac{6^2}{((n\eta)^3-n\eta)^2}\sum_{j=1}^{n\eta} (2j-n\eta-1)^2\\
                &= \underbrace{\frac{6^2}{((n\eta)^3-n\eta)^2}}_{{=} \mathrm{O}(n^{-6})}\underbrace{\frac{1}{3}((n\eta)^3-n\eta)}_{{=} \mathrm{O}(n^3)}{=} \mathrm{O}(n^{-3}) {=} \mathrm{O}(n^{-1}).
			\end{align*}
  		
        \section{Proof of convergence of $\Upsilon^{(n)}$ to $\Upsilon$}\label{appendix:limitprocessupsilon}
		
		\subsection{The limit process}\label{sec:def_proc}
        Recall that we fix $T>0$ and $0<\eta<\frac{T}{2}$, denote $\para:=[\eta,T-\eta]$ and define a centered Gaussian process $\Upsilon=(\Upsilon_t)_{t\in\para}$ by its covariance function $C:\para\times \para\to\Rset$ defined in \eqref{eq:limit_covariance}. By a criterion of Fernique \cite{fer64} we obtain that we have a version of $\Upsilon$ with continuous sample paths. (Choose, e.g., $\psi(x):=6(x/\eta)^3$ in the formulation of Marcus and Shepp \cite{mar_she_70}.) The space of continuous functions on $\para$ is denoted by $C(\para)$.
		\begin{remark}
			We will not make use of the fact that $(\Upsilon_t)_{t\in\para}$ has a representation as an It\^{o}-integral with respect to Brownian motion with a deterministic kernel: With the functions $g_u$ defined in (\ref{def_v_function}) we have, in distribution,
			\begin{align*}
				(\Upsilon_t)_{t\in\para}=\Bigg( \sqrt{3/2}\eta^{-3/2}\int_{t-\eta}^{t+\eta} g_t(v)\,dB_v\Bigg)_{t\in\para},
			\end{align*}
			where $(B_v)_{v\in [0,\infty)}$ denotes standard Brownian motion.
		\end{remark}
		
		\subsection{The discrete processes for fixed $\mathbf{n}$}
        Recall that the processes $\Upsilon^{(n)}$ are defined, see also \eqref{rep_gamma_null}, by $\Upsilon^{(n)}_t := \Gamma_{n\eta,nt}^{(1)}$ with
		\begin{align*}
            \Gamma_{n\eta,nt}^{(1)}&=\left(\frac{2}{3}((n\eta)^3-n\eta)\right)^{-1/2} \\
            &\quad\quad~\times\Bigg(\sum_{j=1}^{n\eta}(2j-n\eta-1)\left(Z_{\lfloor nt\rfloor+j}^{(1)}-Z_{\lfloor nt\rfloor-n\eta+j}^{(1)}\right)\Bigg)
		\end{align*}
        and considered as processes in the c\`{a}dl\`{a}g space $D(\para)$ endowed with the topology induced by the supremum norm $\|\,\cdot\,\|_\infty$. Note that they are measurable. To prepare for the process convergence of the $\Upsilon^{(n)}$ in Proposition \ref{proposition:convergenceupsilonn} we first study the covariance functions of $\Upsilon^{(n)}$ and $\Upsilon$.

		\subsection{The covariance functions}\label{section:covariancefunctionupsilon}
		Recall the definition of $C$ in (\ref{eq:limit_covariance}) and further denote the covariance function of $\Upsilon^{(n)}$ by
		\begin{align*}
			C_n(s,t):=\mathrm{Cov}(\Upsilon^{(n)}_s,\Upsilon^{(n)}_t)
			=\E[\Upsilon^{(n)}_s\Upsilon^{(n)}_t].
		\end{align*}
		We show uniform convergence of the covariance functions as stated in Proposition \ref{lemma:convergencecovariancefunction}.
  
		\noindent{\em Proof of Proposition \ref{lemma:convergencecovariancefunction}.}
		In order to obtain a convenient form of the covariance functions of the $\Upsilon^{(n)}$ we note that, see (\ref{eq:upsiloncontaininggu2}), we have
		\begin{align}
			\Upsilon^{(n)}_t
			&=\sqrt{\frac{3}{2}}\left((n\eta)^3-n\eta\right)^{-1/2}\sum_{j=-\infty}^{\infty} g^{(n)}_t(j)V_j,
			\label{ref_proc}
		\end{align}
		where $\{V_j\,|\,j\in\Zset\}$ is a set of independent standard normal $\mathcal{N}(0,1)$ distributed random variables and, for $u \in\Rset$, functions $g^{(n)}$ are defined by 
		\begin{align}
  g^{(n)}_u: \Zset \to \Rset \label{def_gnu}
  		\end{align}
    with
    		\begin{align*}
			j\mapsto \left\{ \begin{array}{cl}
				2|j-\lfloor nu\rfloor|-n\eta +\mathrm{sgn}(\lfloor nu\rfloor-j), & 0<|\lfloor nu\rfloor -j| < n\eta\\ &\mbox{and } j<\lfloor nu\rfloor \mbox{ or}\\
                &0<|\lfloor nu\rfloor -j|\leq n\eta\\ 
                &\mbox{and } j>\lfloor nu\rfloor,\nonumber\\
                 2|j-\lfloor nu\rfloor|-n\eta +1, &\mbox{for } j = \lfloor nu \rfloor \nonumber, \\
				0, & \mbox{otherwise},\nonumber
			\end{array} \right.
		\end{align*}
		with
		\begin{align*}
			\mathrm{sgn}(x):= \left\{  \begin{array}{cl}
				1, &\mbox{for } x>0, \\
				0, &\mbox{for } x=0, \\
				-1, &\mbox{for } x<0.
			\end{array} \right.
		\end{align*}
		Since  $\{V_j\,|\,j\in\Zset\}$ is a set of independent standard normal random variables we obtain
		\begin{align}
			C_n(s,t)&=\frac{3}{2((n\eta)^3-n\eta)}\sum_{j=-\infty}^{\infty}g^{(n)}_s(j) g^{(n)}_t(j)\nonumber\\
   &= \frac{3}{2\eta^{3}+\mathrm{O}(1/n)}\Bigg\{ \frac{1}{n}\sum_{j=-\infty}^{\infty}\left(\frac{1}{n}g^{(n)}_s(j)\right)\left( \frac{1}{n}g^{(n)}_t(j)\right)\Bigg\}.\label{riemann_summe}
		\end{align}
		We interpret the expression within $\{\,\cdot\,\}$ in (\ref{riemann_summe}) as a Riemann sum. For this, we define functions $g_u$ for $u\in \Rset$ by
				\begin{align}\label{def_v_function}
			g_u:\Rset\to \Rset ,\qquad x\mapsto \left(2|x-u|-\eta\right)\mathbf{1}_{[u-\eta,u+\eta]}(x),
		\end{align}
				where $\mathbf{1}_A$ denotes the indicator function of a set $A$. Now, the expression within $\{\,\cdot\,\}$ in (\ref{riemann_summe}) is a Riemann sum for the integral  $\int_{-\infty}^{\infty}g_s(x)g_t(x)\,dx$ up to summands of the order $\mathrm{O}(1/n)$, using the big-$\mathrm{O}$ notation. We obtain, as $n\to\infty$, that
				\begin{align}\label{riemann_approx}
			\left|\frac{1}{n}\sum_{j=-\infty}^{\infty}\left(\frac{1}{n}g^{(n)}_s(j)\right)\left( \frac{1}{n}g^{(n)}_t(j)\right)- \int_{-\infty}^{\infty}g_s(x)g_t(x)\,dx\right| = \mathrm{O}\left(\frac{1}{n}\right)
		\end{align}
        uniformly within $s,t\in \para$. An evaluation of the integral $\int_{-\infty}^{\infty}g_s(x)g_t(x)\,dx$ can be done by slightly tedious but elementary calculations and yields 		
		\begin{align}\label{eval_int}
			\frac{3}{2\eta^3} \int_{-\infty}^{\infty}g_s(x)g_t(x)\,dx= \kappa\left(\frac{|s-t|}{\eta}\right)=C(s,t)
		\end{align}
        with $\kappa$ and $C$ given in (\ref{eq:covariance_kappa}) and (\ref{eq:limit_covariance}) respectively. Combining (\ref{riemann_summe}), (\ref{riemann_approx}) and (\ref{eval_int}) we obtain as $n\to\infty$ that
		\begin{align*} 
			C_n(s,t)=\mathrm{Cov}(\Upsilon^{(n)}_s,\Upsilon^{(n)}_t) \to \mathrm{Cov}(\Upsilon_s,\Upsilon_t)=C(s,t)
		\end{align*}
				uniformly within $s,t\in \para$. \hfill $\square$
						
				The Gaussian processes $\Upsilon$ and $\Upsilon^{(n)}$ induce canonical metrics on $\para$ by
		\begin{align*}
			d_\infty(s,t)&:=\|\Upsilon_t-\Upsilon_s\|_2
			:=\E[|\Upsilon_t-\Upsilon_s|^2]^{1/2},\\
			d_n(s,t)&:=\|\Upsilon^{(n)}_t-\Upsilon^{(n)}_s\|_2.
		\end{align*}
		Later, we need an upper bound on the $d_n$ for $|s-t|$ small.
		\begin{proposition}\label{bound_cano}
			There exists a constant $0<c_\eta<\infty$, depending on $\eta$, such that for the canonical distances of $\Upsilon$ and $\Upsilon^{(n)}$ we have
			\begin{align*}
				d_n(s,t) \le c_\eta |s-t|^{1/2},
			\end{align*}
			for all $s,t\in\para$ with $|s-t|\le \eta$ and all $n\in\Nset \cup\{\infty\}$.
		\end{proposition}
     \noindent   {\em Proof.} For $s,t\in\para$ with $|s-t|\le \eta$ we obtain with (\ref{eq:covariance_kappa}) and (\ref{eq:limit_covariance}) that
		\begin{align}
			d_\infty(s,t)&=\left(\Var(\Upsilon_t-\Upsilon_s)\right)^{1/2}=\left(2(1-C(s,t))\right)^{1/2}\nonumber\\
            &=\left(2\left(\frac{3}{2}\frac{|s-t|}{\eta}+ 3\frac{|s-t|^2}{\eta^2}-3\frac{|s-t|^3}{\eta^3}\right)\right)^{1/2}\label{con_metr}\nonumber\\
			&\le \left(3\frac{|s-t|}{\eta}+ 6\frac{|s-t|}{\eta}\right)^{1/2}=\frac{3}{\sqrt{\eta}} |s-t|^{1/2}.
		\end{align}	
		Note that $\Var(\Upsilon_t)=\Var(\Upsilon^{(n)}_t)=1$ for all $t$ implies, using notation introduced in the proof of Proposition \ref{lemma:convergencecovariancefunction}, that
		\begin{align*}
			\frac{1}{n}\sum_{j=-\infty}^{\infty}\left(\frac{1}{n}g^{(n)}_t(j)\right)\left( \frac{1}{n}g^{(n)}_t(j)\right)= \int_{-\infty}^{\infty}g_t(x)g_t(x)\,dx=1.
		\end{align*}
		Hence, by triangle inequality, (\ref{riemann_summe}) and (\ref{riemann_approx}), we obtain
		\begin{align*}
			\lefteqn{\frac{2\eta^3+\mathrm{O}(1/n)}{3}|C_n(s,t)-C(s,t)|} \nonumber\\
			&\le \left|\frac{1}{n}\sum_{j=-\infty}^{\infty}\left(\frac{1}{n}g^{(n)}_s(j)\right)\left( \frac{1}{n}g^{(n)}_t(j)\right)- \frac{1}{n}\sum_{j=-\infty}^{\infty}\left(\frac{1}{n}g^{(n)}_t(j)\right)\left( \frac{1}{n}g^{(n)}_t(j)\right)\right|\\
   &\quad~+ \left|\int_{-\infty}^{\infty}g_t(x)g_t(x)\,dx - \int_{-\infty}^{\infty}g_s(x)g_t(x)\,dx
			\right| \nonumber \\
			&\le \frac{1}{n}\sum_{j=-\infty}^{\infty}\left|\frac{1}{n}g^{(n)}_s(j)-\frac{1}{n}g^{(n)}_t(j)\right| \frac{1}{n}|g^{(n)}_t(j)| +\int_{-\infty}^{\infty}|g_t(x)||g_t(x)-g_s(x)|\,dx \nonumber\\
			&\le c|s-t|,
		\end{align*}
		with appropriate $0<c<\infty$ uniformly in $n$, noting that $|\frac{1}{n}g^{(n)}_s(j)-\frac{1}{n}g^{(n)}_t(j)|\le 2|s-t|$ for all $j=\lfloor n(s\wedge t)\rfloor - n\eta,\ldots,\lfloor n(s\vee   t)\rfloor + n\eta$ and $|g_t(x)-g_s(x)|\le 2|s-t|$ for all $x\in [(s\vee t) -\eta, (s\wedge t) +\eta]$, and using trivial bounds outside these summation and integration ranges, respectively. Hence, for $s,t\in\para$ we have
		\begin{align}\label{appr_zero}
			|C_n(s,t)-C(s,t)|\le \frac{3c}{2\eta^3+\mathrm{O}(1/n)}|s-t|.
		\end{align}
		Using (\ref{appr_zero}) in (\ref{con_metr}) we further obtain for all $|s-t|\le \eta$ that
		\begin{align*}
			\lefteqn{d_n^2(s,t)}\\
   &=2(1-C_n(s,t))\\
   &\le 2\Bigg(\frac{3}{2}\frac{|s-t|}{\eta}+ 3\frac{|s-t|^2}{\eta^2}-3\frac{|s-t|^3}{\eta^3}+\frac{3c}{2\eta^3+\mathrm{O}(1/n)}|s-t|)\Bigg)\\
   &\le c_\eta^2 |s-t|
		\end{align*}
		for an appropriate $c_\eta\ge 3/\sqrt{\eta}$. \hfill $\square$

		\subsection{Convergence of the processes} \label{app_func_proof}
		We now show process convergence as stated in Proposition \ref{proposition:convergenceupsilonn}. To prove Proposition \ref{proposition:convergenceupsilonn} our bounds on the covariance functions from Section \ref{section:covariancefunctionupsilon} allow to verify conditions of general theorems which we first restate in a notation and form suitable for further use and tailored to our setting. Firstly, we use a criterion from Pollard \cite[Sec V, Thm 3]{Pollard1984}:
		\begin{theorem}\label{pollard_crit}
			Let $\kgs{\Xi},\kgs{\Xi}^{(n)}$, $n\ge 1$, be random elements in $(D(\para),\|\,\cdot\,\|_\infty)$. Suppose that $\Prob(\kgs{\Xi}\in C(\para))=1$. Necessary and sufficient for   $\kgs{\Xi}^{(n)} \stackrel{d}{\longrightarrow} \kgs{\Xi}$ in $(D(\para),\|\,\cdot\,\|_\infty)$ as $n\to\infty$, is that the following two conditions hold:
			\begin{enumerate}
				\item[(i)]
				We have $\kgs{\Xi}^{(n)}\stackrel{\mathrm{fdd}}{\longrightarrow} \kgs{\Xi}$ as $n\to \infty$,
				\item[(ii)]
				For all $\varepsilon, \delta>0$ there exist $\eta=t_0<t_1<\cdots<t_m=T-\eta$ such that
				\begin{align}\label{straff}
					\limsup_{n\to\infty} \Prob\left(\max_{i=0,\ldots,m-1}\sup_{t\in[t_i,t_{i+1})}|\kgs{\Xi}^{(n)}_t-\kgs{\Xi}^{(n)}_{t_i}|>\delta\right)<\varepsilon.
				\end{align}
			\end{enumerate}
		\end{theorem}
		
		Secondly, in order to verify condition (ii) in  Theorem \ref{pollard_crit} we apply a bound of Dirksen \cite[Theorem 3.2]{dir15} which is based on Talagrand's \cite{ta05,ta21} generic chaining method. Dirksen showed that there exist constants $0<C,D<\infty$ such that for any real-valued centered Gaussian process $(\kgs{\Xi}_t)_{t\in \tau}$ with appropriate metric space $(\tau,d)$ and any $1\le p <\infty$ and $t_0\in\tau$ we have
		\begin{align}\label{dirk}
			\E\left[\sup_{t\in\tau}\left|\kgs{\Xi}_t-\kgs{\Xi}_{t_0}\right|^p\right]^{1/p}\le  C\gamma_{2,p}(\tau,d) + 2\sup_{t\in\tau}\E\left[\left|\kgs{\Xi}_t-\kgs{\Xi}_{t_0}\right|^p\right]^{1/p}.
		\end{align}
		Here, $d(s,t)=\E[|\kgs{\Xi}_s-\kgs{\Xi}_t|^2]^{1/2}$ is the canonical metric of $(\kgs{\Xi}_t)_{t\in \tau}$ and  $\gamma_{2,p}(\tau,d)$ is recalled below. The bound (\ref{dirk}) applies to our processes $\Upsilon, \Upsilon^{(n)}$, i.e., we have appropriate $(\tau,d)$, since our processes have continuous respectively piecewise constant sample paths, cf.~Remark 3.1 in \cite{dir15}.
		
		Based on these results we obtain Proposition \ref{proposition:convergenceupsilonn} as follows:
		
		\noindent{\em Proof of Proposition \ref{proposition:convergenceupsilonn}.}
		It is sufficient to verify the conditions of Theorem \ref{pollard_crit} with $\kgs{\Xi}^{(n)}=\Upsilon^{(n)}$ and $\kgs{\Xi}=\Upsilon$. First note that we have $\Prob(\Upsilon\in C(\para))=1$ as noted in Section \ref{sec:def_proc}.

        Condition (i) in Theorem \ref{pollard_crit} follows from Proposition \ref{lemma:convergencecovariancefunction}, as multivariate random Gaussian vectors with converging covariance matrices converge in distribution to the random Gaussian vector with limiting covariance matrix, e.g., by L\'evy's continuity theorem.
		
        To verify condition (ii) in Theorem \ref{pollard_crit} let $\varepsilon, \delta>0$. We choose $t_j:= \eta+(T-2\eta)j/m$ for $j=0,\ldots,m$ with $m$ to be determined later. Note that by subadditivity, Markov's inequality and that $\Upsilon^{(n)}$ is stationary we have
		\begin{align*}
			\lefteqn{\Prob\left(\max_{i=0,\ldots,m-1}\sup_{t\in[t_i,t_{i+1})}|\Upsilon^{(n)}_t-\Upsilon^{(n)}_{t_i}|>\delta\right)}\\
   &\le m \Prob\left(\sup_{t\in[0,T/m]}|\Upsilon^{(n)}_t-\Upsilon^{(n)}_{0}|>\delta\right)\\
            &\le \frac{m}{\delta^4}\E\left[\sup_{t\in[0,T/m]}\left|\Upsilon^{(n)}_t-\Upsilon^{(n)}_{0}\right|^4  \right].
		\end{align*}
		Hence, for (\ref{straff}) to be satisfied for the processes $\Upsilon^{(n)}$ in view of the latter display it is sufficient to show that
		\begin{align}\label{straff_b}
			\sup_{n\ge 1} \E\left[\sup_{t\in[0,T/m]}\left|\Upsilon^{(n)}_t-\Upsilon^{(n)}_{0}\right|^4\right]= \mathrm{O}\left(\frac{1}{m^2}\right)\quad (m\to\infty).
		\end{align}
			In order to obtain (\ref{straff_b}) we bound the right hand side of (\ref{dirk}) with $\tau=\tau_m:=[0,T/m]$, $p=4$ and $d=d_n$ the canonical distance of $\Upsilon^{(n)}$. To bound the second summand on the right hand side of (\ref{dirk}) recall that the fourth moment (kurtosis) of a standard normal random variable has the value $3$. For $t_0,t\in \tau_m$ and $m>T/\eta$, using Proposition \ref{bound_cano}, we find that $\Upsilon^{(n)}_t-\Upsilon^{(n)}_{t_0}$ has a centered normal distribution with variance $d_n^2(t,t_0)\le c_\eta^2 |t-t_0|\le c_\eta^2 T/m$. Hence, the fourth moment of $\Upsilon^{(n)}_t-\Upsilon^{(n)}_{t_0}$ is upper bounded by $3c_\eta^4 T^2/m^2$ and we obtain for all $m>T/\eta$ that
		\begin{align}\label{dirk_app1}
			2\sup_{t\in\tau_m}\E\left[\left|\Upsilon^{(n)}_t-\Upsilon^{(n)}_{0}\right|^4\right]^{1/4}\le \frac{2\sqrt[4]{3T^2}c_\eta}{\sqrt{m}}=\mathrm{O}\left(\frac{1}{\sqrt{m}}\right),
		\end{align}
		where the big-$\mathrm{O}$-constant is independent of $n$.
		
		To bound the first summand on the right hand side of (\ref{dirk}) we recall the definition of $\gamma_{2,p}(\tau,d)$ from \cite{dir15}: A sequence $\mathcal{T}=(T_\ell)_{\ell\ge 0}$ of subsets $T_\ell\subset \tau$ are called admissible if their cardinalities satisfy $|T_0| =1$ and $|T_\ell|\le 2^{2^\ell}$ for all $\ell\ge 1$. Then, with $p=4$, we have
			\begin{align*}
			\gamma_{2,4}(\tau,d) := \inf_\mathcal{T}\sup_{t\in \tau}\sum_{\ell =2}^\infty 2^{\ell/2} d(t,T_\ell),
		\end{align*}
	    where the infimum is taken over all admissible sequences $\mathcal{T}$. With  $\tau=\tau_m=[0,T/m]$ and $d=d_n$ we simply define $T_0:=\{0\}$ and
		\begin{align*}
			T_\ell:=\{(jT)/(m2^{2^\ell}):j=1\ldots,2^{2^\ell}\}, \quad \ell\ge 1.
		\end{align*}
		This yields an admissible sequence $(T_\ell)_{\ell\ge 0}$. For all $t\in \tau_m$ and $\ell\ge 1$ there exists $s\in T_\ell$ with $|s-t|\le T/ (m2^{2^\ell})$. For all $m$ sufficiently large so that $T/m\le \eta$ we obtain from Proposition \ref{bound_cano} that
		\begin{align*}
			d_n(t,T_\ell)\le d_n(t,s)\le c_\eta (T/ (m2^{2^\ell}))^{1/2}= c_\eta\sqrt{T/m}2^{-2^{\ell-1}}.
		\end{align*}
		Hence, for all $m\ge T/\eta$ we have
		\begin{align}\label{dirk_app2}
			\gamma_{2,4}(\tau_m,d_n)&\le \sum_{\ell=2}^\infty 2^{\ell/2}c_\eta\sqrt{T/m}2^{-2^{\ell-1}}\nonumber\\
   &\le \frac{c_\eta\sqrt{T}}{\sqrt{m}}\sum_{\ell=2}^\infty 2^{-2^{\ell-1}+\ell/2}=\mathrm{O}\left(\frac{1}{\sqrt{m}}\right),
		\end{align}
        since the latter series is convergent. Note, that the big-$\mathrm{O}$-bound is independent of $n$. Combining (\ref{dirk}), (\ref{dirk_app1}) and (\ref{dirk_app2}) implies (\ref{straff_b}) and hence  condition (ii) in Theorem \ref{pollard_crit}. \hfill $\square$

	\end{appendix}

\subsubsection*{Author contributions statement}
     S.P. writing -- original draft (equal), analysis (lead), data curation (equal), theoretical proofs (equal), visualization (equal), T.E.: Recording of data set (equal), data curation (equal), writing -- original draft (supporting) P.G.: Recording of data set (equal), data curation (supporting) E.S.: supervision (supporting), funding acquisition (lead), writing -- original draft (supporting) R.N.: theoretical proofs (lead), writing -- original draft (equal) G.S.: writing -- original draft (equal), theoretical proofs (equal), formal analysis (supporting), funding accquisition (lead), visualization (equal), supervision (lead)



\subsubsection*{Acknowledgement}
We thank Jaclyn-Katrin Boppenmaier for investigation of the behavior of the leaf plot in various settings and two anonymous reviewers for their constructive comments.

\subsubsection*{Funding}
This work was supported by the LOEWE Schwerpunkt CMMS - Multi-scale Modelling in the Life Sciences.


\subsubsection*{Supplements}
Supplement 1: {\em Univariate change point detection:} An elaboration on the univariate approach on change point detection in the movement direction $\vartheta$ in an LW with constant step length $r$. The supplement includes simulations on the significance level and testpower and a discussion on limitations in cases where the step length shows change points.\\

\noindent
Supplement 2: {\em Unbiasedness of MLEs:} Proof of the Unbiasedness of the MLEs from Proposition \ref{theorem:ml_mu_positions}.\\

\noindent
Supplement 3: {\em Proof of Strong Consistency of MLEs:} An adapted SLLN and proof of strong consistency of $\hat b$ and $\hat \sigma^2$ from Proposition 3.1.\\

\noindent
Supplement 4: {\em R-Code:} Not contained in the arXiv-version. Please contact the corresponding author Gaby Schneider: \href{mailto:schneider@math.uni-frankfurt.de}{schneider@math.uni-frankfurt.de}



\bibliographystyle{plain} 

\section*{Supplement 1: Univariate change point detection}
	
	Here we shortly discuss problems and restrictions of the univariate moving window approach. As described there, we derive the process $(D^{\kgs{\vartheta}}_i)_i$ of differences of estimated directions according {to the following equation}
    \begin{equation}
        D_i^{\kgs{\vartheta}}:=D_i(\hat{\vartheta}_{\ell}\kgs{(i)},\hat{\vartheta}_{r}\kgs{(i)}):=
		\text{atan2}(\sin(\hat{\vartheta}_{r}\kgs{(i)}-\hat{\vartheta}_{\ell}\kgs{(i)}),\cos(\hat{\vartheta}_{r}\kgs{(i)}-\hat{\vartheta}_{\ell}\kgs{(i)})),
	\end{equation}
    and use its maximum, $M$, as a test statistic for the null hypothesis of no change in direction. In order to derive a rejection threshold for the statistical test, we need the distribution of $M$ under the null hypothesis. As it is not easily accessible formally, we simulate the distribution of $M$ by assuming that all other parameters, particularly $r$, are constant. The empirical $(1-\alpha)$-quantile of the resulting simulations then serves as rejection threshold, $Q$.
	
\begin{figure}[!h]
		\centering
		\includegraphics[width=1\textwidth]{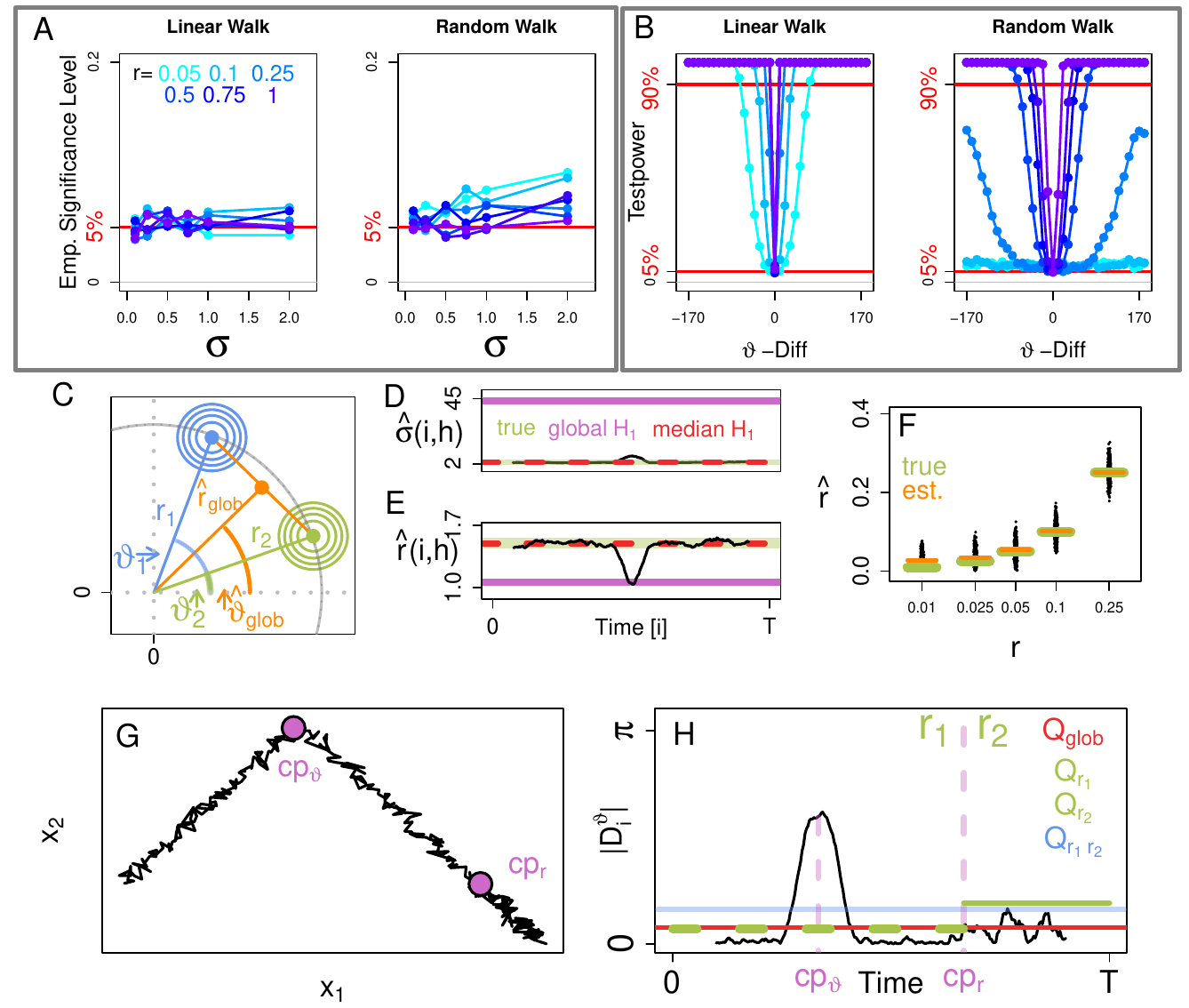}
 		\caption{Empirical significance level (A) and test power (B) for the LW and RW for $\alpha=5\%$. $1000$ simulations were performed with $h=30, \vartheta=35\kgs{/180\cdot\pi}$ and $T=400$, and $\sigma=0.5$ for the testpower. (C) Distribution of increments before (blue) and after (green) the change point. Estimated expectation from mixture of increments (orange) and the resulting estimated direction and step length. Local and global estimators of the model parameters step length $r$ (D) and $\sigma^2$ (E) with window size $h=30$ from a track with one change in direction. (F) Estimation of the step length $r$ from data without change points. Mean of estimated step length (green line) and true step length (orange line). 1000 Simulations per step length $r$, $\sigma^2=1$, $h=30.$ (G) Track from the Linear Walk model with one change point $cp_\vartheta$ in direction and one change point $cp_r$ in step length. (H) Process of $|D_t|$ of the absolute value of differences of directions (black) derived from the track in G, threshold $Q_{glob}$ derived from global estimation of step length (red). Thresholds $Q_{r_1}$ and $Q_{r_2}$ derived from estimations of step lengths on either side of the change point in step length (green). Threshold $Q_{r_1,r_2}$ derived from tracks with the change point in step length (blue).}	
		\label{fig:uni5}
	\end{figure}

	Figure \ref{fig:uni5} A shows that for the used parameter sets, the chosen significance level can be kept in  this setting both for the RW and the LW already for a window size of $h=30$. However, one should note that for small sample sizes, the step length $r$ may tend to be overestimated. To see this, consider first the degenerate case with $r=0$ in which the distribution of increments is centered at the origin. The distribution of $\hat \mu$ is then also normal and centered at the origin. Thus, the true expectation of $r$ should be zero, but by construction, all estimates $\hat r$ are positive, such that $r$ must necessarily be overestimated. If the distribution of $\hat \mu$ is shifted only slightly away from zero, this effect will remain. It is therefore necessary to assume a sufficiently large value of $\frac{r}{\sigma}$ (Figure \ref{fig:uni5} F).

	If we apply the global parameter estimates in cases with change points in the model parameters, we observe a certain bias both in the RW and in the LW case. If a track contains one change point in direction, the increments between successive points are assumed to originate from two different normal distributions with the same variance $\sigma^2$ and the same distance $r$ to the origin, but with different directions from the origin (Figure \ref{fig:uni5} C). In that case, the step length $r$ will tend to be underestimated (purple line in panel E), and the standard deviation $\sigma$ will tend to be overestimated  (purple line in panel D). 	We therefore propose to apply a moving window approach and to use the median of the resulting process as a parameter estimate 		(red dashed line in Figure \ref{fig:uni5} D and E).

	In our simulations with change points, the test power increased with the difference in directions and with increasing step length. In addition, it was higher for the LW than for the RW (Figure \ref{fig:uni5} B) for the same parameters, which is due to a higher variability in the RW.
		
	This indicates that the proposed statistical tests are suitable for detecting changes in the movement direction under the assumption of constant step length. However, one should note that additional change points in $r$ as illustrated in Figure \ref{fig:uni5} G will affect the rejection threshold $Q$, because $Q$ depends on $r$ (Figure \ref{fig:uni5} H). Smaller values of $r$ yield a higher variability of $|D|$ and thus, a higher $Q$ (green lines in Figure \ref{fig:uni5} H). Ignoring potential changes in $r$, one will obtain a false global estimate of  $r$, such that the resulting global rejection threshold $Q_{glob}$ (red line in Figure \ref{fig:uni5} H) cannot adhere to the chosen significance level.

\section*{Supplement 2: Unbiasedness of MLE based estimates from Proposition \textcolor{black}{2.3}}
	
	Here we show the unbiasedness of the estimators derived in Proposition \textcolor{black}{2.3}, which states that for a uni-directional \textbf{LW}($\vartheta$,$r$,$b$,$\sigma$,$\emptyset$) with a sequence of positions $X_{i+1}, X_{i+2}, \ldots, X_{i+h}$, the MLEs for $\mu$ and $\sigma^2$ are given, respectively, by 
		\begin{align*}
			\hat{\mu}(i,h)&=\frac{6}{h^3-h}\sum_{j=1}^{\lfloor\frac{h}{2}\rfloor}\left(h-(2j-1)\right)\left(X_{i+j}-X_{i+h-j+1}\right), \\
			\hat{b}(i,h)  &= \bar X\kgs{_i}- (t+(h+1)/2)\hat\mu, \quad \text{and} \\
			\hat{\sigma}^{2}(i,h) &=\frac{1}{2h}\sum_{j=1}^{h} \| X_{i+j}-(i+j)\hat\mu-\hat b\|^2,
			\end{align*}
	where $\bar{X}\kgs{_i}:=h^{-1}\sum_{j=1}^h{X_{i+j}}$.
	
	We first note the following equalities for the weights in $\hat{\mu}^{(d)}$ and $\hat{b}^{(d)}$, which we denote by $w_\mu(h,j):=2j-h-1$ and $w_b(h,i,j):=-6hj-12ji-6j+4h^2+6hi+6h+6i+2$, respectively. Then we have
\begin{align}
		&\sum_{j=1}^{h}w_\mu(h,j)=0 \quad\text{ and }\quad \sum_{j=1}^{h}w_\mu(h,j)(i+j)=\frac{h^3-h}{6}, \label{eq:weightsumsmu}\\
		&\sum_{j=1}^{h}w_b(h,i,j)=h^3-h\quad \text{ and }\quad \sum_{j=1}^{h}w_b(h,i,j)(i+j)=0. \label{eq:weightsumsb}
	\end{align}
	
	We then find for $\hat{\mu}^{(d)}$
\begin{align*}
		\mathbb{E}\left[\hat{\mu}^{(d)}\right]&= \frac{6}{h^3-h}\sum_{j=1}^{h}w_\mu(h,j)\mathbb{E}\left[X_{i+j}^{(d)}\right]\\
		&=\frac{6}{h^3-h}\sum_{j=1}^{h}w_\mu(h,j)\left((i+j)\mu^{(d)}+b^{(d)}\right)\\
		&=\frac{6}{h^3-h}\left(\mu^{(d)}\sum_{j=1}^{h}w_\mu(h,j)(i+j) +b^{(d)}\sum_{j=1}^{h}w_\mu(h,j)\right)\\
		&=\frac{6}{h^3-h}\frac{h^3-h}{6}\mu^{(d)} = \mu^{(d)},
	\end{align*}
	where we use equation (\ref{eq:weightsumsmu}). For $\hat{b}^{(d)}$  we use equation (\ref{eq:weightsumsb}) and observe
\begin{align}
		\mathbb{E}\left[\hat{b}^{(d)}\right]&= \nonumber \frac{1}{h^3-h}\sum_{j=1}^{n}w_b(h,i,j)\mathbb{E}\left[X_{i+j}^{(d)}\right]\nonumber\\
		&= \frac{1}{h^3-h}\sum_{j=1}^{n}w_b(h,i,j) \left((i+j)\mu^{(d)}+b^{(d)}\right)\nonumber \\\
		&= \frac{1}{h^3-h}\left(\mu^{(d)}\sum_{j=1}^{n}w_b(h,i,j)(i+j)+b^{(d)}\sum_{j=1}^{n}w_b(h,i,j)\right) = b^{(d)}. \label{eq:b_unbiased}
	\end{align}
    For the expectation of $\textcolor{blue}{\hat\sigma}^2$ we first insert equation \eqref{eq:MLEbLM} into equation \eqref{eq:MLEsigmaLM}, noting that
    \begin{align*}
        X_{i+j} -(i+j)\hat{\mu}-\hat{b} &= X_{i+j} - (i+j)\hat{\mu} - \bar{X} +\left(i+\frac{h+1}{2}\right)\hat{\mu}\\
        &= X_{i+j} - \bar{X} + \hat{\mu}\left( \frac{h+1}{2}-j \right).
    \end{align*}
    Now using $X_{i+j}=(i+j)\mu+b+\sigma Z_{i+j}$, we note that 
    \begin{align*}
        X_{i+j}-\bar{X} = \mu\left( j-\frac{h+1}{2}\right)+ \sigma\left(Z_{i+j}-\bar{Z} \right),
    \end{align*}
    where $\bar{Z}=\frac{1}{h}\sum_{j=1}^hZ_{j+i}$. For unbiasedness of $\hat{\sigma}^2$ we thus need to show that 
    \begin{align}
        \mathbb{E}\left[\sum_{j=1}^h \left\lVert\sigma\left( Z_{i+j}-\bar{Z}\right)+\left(\hat{\mu}-\mu\right)\left( \frac{h+1}{2}-j\right) \right\rVert^2  \right]= (2h-4)\sigma^2.\label{eq:hatsigmasquared_alternativeform}
    \end{align}
    The expression on the left of \eqref{eq:hatsigmasquared_alternativeform} equals $\sum_{j=1}^h\sum_{d=1,2}(a_1+a_2+a_3)$, with
    \begin{align*}
    a_1&= \sigma^2\mathbb{E}\left[\left( Z_{i+j}^{(d)}-\bar{Z}^{(d)}\right)^2\right],\\
    a_2&= 2\sigma\left( \frac{h+1}{2}-j\right)\mathbb{E}\left[\left(Z_{i+j}^{(d)}-\bar{Z}^{(d)}\right)\left(\hat{\mu}^{(d)}-\mu^{(d)}\right) \right],\\
    a_3&=\left( \frac{h+1}{2}-j\right)^2\mathbb{E}\left[\left(\hat{\mu}^{(d)}-\mu^{(d)}\right)^2\right].\label{eq:unbiasednesshatsigmasquared}
    \end{align*}
    We then observe
    \begin{align*}
        a_1&= \sigma^2\mathbb{E}\left[Z_{i+j}^{(d)^2}-2Z_{i+j}^{(d)}\bar{Z}^{(d)}+\bar{Z}^{(d)^2} \right]=\sigma^2\left(1-\frac{1}{h}\right),\\
        a_2 &= 2\sigma\left(\frac{h+1}{2}-j\right)\left(\mathbb{E}\left[Z_{i+j}^{(d)}\hat{\mu}^{(d)}\right] - \mathbb{E}\bigg[\bar{Z}^{(d)}\hat{\mu}^{(d)}\bigg] \right).
     \end{align*}
Now we note that
    \begin{align*}
        \mathbb{E}\left[Z_{i+j}^{(d)}\hat{\mu}^{(d)}\right]&= \mathbb{E}\left[Z_{i+j}^{(d)}\frac{6}{h^3-h}\sum_{k=1}^hw_\mu(h,k)X_{i+k}^{(d)}  \right]\\
        &=  \frac{6}{h^3-h}w_\mu(h,j)\mathbb{E}\left[Z_{i+j}^{(d)}X_{i+j}^{(d)}\right] \\
        &= \frac{6}{h^3-h}w_\mu(h,j)\mathbb{E}\left[Z_{i+j}^{(d)}(\mu^{(d)}+b^{(d)}+\sigma Z_{i+j}^{(d)})\right]\\
        &= \sigma \frac{6}{h^3-h}w_\mu(h,j)
    \end{align*}
    and
    \begin{align*}  \mathbb{E}\left[\bar{Z}^{(d)}\hat{\mu}^{(d)}\right] &= \mathbb{E}\left[\frac{1}{h}\sum_{k=1}^hZ_{i+k}^{(d)}\frac{6}{h^3-h}\sum_{l=1}^hw_\mu(h,l)X_{i+l}^{(d)}\right]\\
        &=\frac{6}{h(h^3-h)} \mathbb{E}\left[\sum_{k=1}^h\sum_{l=1}^hw_\mu(h,l)X_{i+l}^{(d)}Z_{i+k}^{(d)} \right]\\
        &= \frac{6}{h(h^3-h)} \left( \sum_{l=1}^h{w_\mu(h,l)\mathbb{E}\left[Z^{\textcolor{blue}{(d)}}_{i+l}X^{\textcolor{blue}{(d)}}_{i+l}\right]}  \right)\\
        &= \frac{6\sigma}{h(h^3-h)}{\sum_{l=1}^hw_\mu(h,l)}= 0
    \end{align*}
and hence find 
\begin{align*}
       a_2 &= 2\sigma\left(\frac{h+1}{2}-j\right)\left( \sigma\frac{6}{h^3-h}w_\mu(h,j) \right)
        =-\sigma^2\frac{6}{h^3-h}w_\mu(h,j)^2.
\end{align*}
For $a_3$ we have
    \begin{align*}
        a_3 
        &= \left(\frac{h+1}{2}-j\right)^2\left(\mathbb{E}\bigg[\hat{\mu}^{(d)^2}\bigg] -2\mu^{(d)}\mathbb{E}\bigg[\hat{\mu}^{(d)}\bigg] + \mu^{(d)^2}\right)
        =\frac{1}{4}\frac{12}{h^3-h}\sigma^2w_\mu(h,j)^2.
    \end{align*}
In total, the left hand side of \eqref{eq:hatsigmasquared_alternativeform} equals
    \begin{align*}
        \lefteqn{\sum_{j=1}^h\sum_{d=1,2}(a_1+a_2+a_3)}\\
        &=\sum_{j=1}^h\sum_{d=1,2} \left(\sigma^2\left(1-\frac{1}{h}\right) -\sigma^2\frac{6}{h^3-h}w_\mu(h,j)^2 + \frac{1}{4}\frac{12}{h^3-h}\sigma^2w_\mu(h,j)^2  \right)\\
        &= 2h\sigma^2\left(1-\frac{1}{h}\right)-2\sigma^2\frac{6}{h^3-h}{\sum_{j=1}^hw_\mu(h,j)^2}+2\sigma^2\frac{3}{h^3-h}\sum_{j=1}^hw_\mu(h,j)^2\\
        &=\left( 2(h-1)-\frac{12}{h^3-h}\frac{1}{3}(h^3-h)+\frac{6}{h^3-h}\frac{1}{3}(h^3-h) \right)\sigma^2\\
        &= (2h-4)\sigma^2,
    \end{align*}
    where we used that $\sum_{j=1}^hw_\mu(h,j)^2=(h^3-h)/3$.
	\qed

\section*{Supplement 3: Strong Consistency of $\hat{b}^{(d)}$ and $\hat{\sigma}^2$}\label{supplement:strongconsistencysigmasquared}

The strong consistency of $\hat{b}^{(d)}$ can be shown using the SLLN for weighted sums from Proposition \textcolor{black}{B.1}. In the time continuous setting, $\hat{b}$ is given by
	\begin{align*}
		\hat{b}(nt,n\eta)&= \frac{1}{(n\eta)^3-n\eta}\sum_{j=1}^{n\eta}w_b(n\eta,nt,j)X^{(d)}_{\lfloor nt\rfloor+j}\\
		&= \sum_{j=1}^{n\eta}\frac{w_b(n\eta,nt,j)}{(n\eta)^3-n\eta}X^{(d)}_{\lfloor nt\rfloor+j}\\
        &= \sum_{j=1}^{n\eta}v(n,j)X^{(d)}_{\lfloor nt\rfloor+i}
	\end{align*}
	where we define 
    \begin{align*}      {w_b(n\eta,nt,j):=}-6n\eta j-12j\lfloor nt\rfloor-6j+4(n\eta)^2+6n\eta \lfloor nt\rfloor+6n\eta+6\lfloor nt\rfloor+2
    \end{align*} and 
$v(n,j):=w_b(n\eta,nt,j)/((n\eta)^3-n\eta)$.
    Again, using the arguments from Appendix \textcolor{black}{B}, it is sufficient to show that $\sum_{j=1}^{m(n)}v(n,j)^2 = \mathrm{O}(n^{-1})$. This is the case because $w_b(n\eta,nt,j) {=} \mathrm{O}(n^2)$ and thus,
	\begin{align*}
		\sum_{j=1}^{n\eta}v(n,j)^2 &= \underbrace{\frac{1}{((n\eta)^3-n\eta)^2}}_{{=} \mathrm{O}(n^{-6})}\underbrace{\sum_{j=1}^{n\eta}\underbrace{w_b(n\eta,nt,j)^2}_{{=}\mathrm{O}(n^4)}}_{{=}\mathrm{O}(n^{5})} {=}\mathrm{O}(n^{-1}).
	\end{align*}       
    To prove strong consistency of $\hat{\sigma}^2$ we use the following variation of the SLLN in a setting of dependent, weighted sums.
    \begin{proposition}[SLLN with weights]\label{proposition:SLLNnew_supp}\quad\\
		Let $X_{jk}$ be centred, i.e., $\mathbb{E}[X_{jk}]=0$, but not necessarily independent or identically distributed random variables with $\sup_{j\ge 1}\mathbb{E}[X_{jk}^6] < \infty$. Let \begin{equation*}S_n:= \sum_{j,k=1}^{m_n} v(n,j,k) X_{jk},\end{equation*} 
where $m_n\in\mathbb{N}$ and $v(n,j,k)\in\mathbb{R}$ are weights. Then $\lim\limits_{n\rightarrow\infty}S_n=0$ almost surely if the following term is in $\mathrm{O}(n^{-2})$
        \begin{align}
            &\sum_{i,j=1}^{m_n}v(n,i,j)^4 + \sum_{\substack{i,j,k,l=1\\ (i,j)\neq(k,l)}}^{m_n} v(n,i,j)^3v(n,k,l)\nonumber\\
            &+\sum_{\substack{i,j,k,l=1\\ (i,j)\neq(k,l)}}^{m_n}v(n,i,j)^2v(n,k,l)^2 +\nonumber\\
            &+\sum_{\substack{i,j,k,l,o,p=1\\ (i,j)\neq(k,l)\\ (i,j)\neq(o,p)\\ (o,p)\neq(k,l)}}^{m_n}v(n,i,j)v(n,k,l)v(n,o,p)^2\nonumber\\
            &+\sum_{\substack{i,j,k,l,o,p,r,s=1;\\ (i,j),(k,l),\\(o,p),(r,s)\\ \text{ pairwise different}}}^{m_n}v(n,i,j)v(n,k,l)v(n,o,p)v(n,r,s).\label{eq:langerterm}
        \end{align}
	\end{proposition}
Proposition \ref{proposition:SLLNnew_supp} can be proven by a standard argument based on the Lemma of Borel--Cantelli and Markov's inequality.\\

Before we show strong consistency of $\hat{\sigma}^2$, we rewrite the estimator. For convenience, these calculations use the time discrete MLE for $\sigma^{2}$ as given in Proposition \textcolor{black}{2.2}. We first refer to a simpler notation of $\hat \sigma^2$ given in \textcolor{black}{(7) in Supplement \kgs{2}}, i.e., 
    \begin{align}
        (2h-4)\hat\sigma^2 &=\sum_{j=1}^h \left\lVert\sigma\left( Z_{i+j}-\bar{Z}\right)+\left(\hat{\mu}-\mu\right)\left( \frac{h+1}{2}-j\right) \right\rVert^2
		\nonumber\\
  &=\sum_{j=1}^h \sum_{d=1,2}\left(\sigma\left( Z_{i+j}^{(d)}-\bar{Z}^{(d)}\right)+\left(\hat{\mu}^{(d)}-\mu^{(d)}\right)\left( \frac{h+1}{2}-j\right)\right)^2, \label{eq:hatsigmasquared_appendix}
	\end{align}
    where the superscript $(d)$ indicates the component of each vector in dimension $d$. For the second summand in the square we find
    \begin{align}
        \lefteqn{\left(\hat{\mu^{(d)}}-\mu^{(d)}\right)\left( \frac{h+1}{2}-j\right)}\nonumber\\ 
        &=\left(\frac{6}{h^3-h}\sum_{k=1}^hw_\mu(h,k)\left((i+k)\mu^{(d)}+b^{(d)}+\sigma Z_{i+k}^{(d)}\right)-\mu^{(d)}\right)\left(-\frac{1}{2}w_\mu(h,j)\right),\label{eq:54}
    \end{align}        
    where $w_{\mu}(h,j):=2j-h-1$. Using that $\sum_{k=1}^hw_\mu(h,k)(i+k)={(h^3-h)}/{6}$ and the notation 
    \begin{align*}
        \hat{\mu}^{(d)}(Z) = \frac{6}{h^3-h}\sum_{k=1}^{h}(2k-h-1)Z^{(d)}_{i+k},
    \end{align*} 
    the term in \eqref{eq:54} equals $-\frac{1}{2}w_\mu(h,j)\sigma \hat{\mu}^{(d)}(Z)$. Thus, we can rewrite $\hat \sigma^2$ as
	\begin{align}
		\hat{\sigma}^{2}(i,h)=\frac{\sigma^2}{2h-4}\sum_{j=1}^{h}\sum_{d=1,2}\left(Z^{(d)}_{i+j}-\bar{Z}^{(d)}-\frac{1}{2}w_{\mu}(h,j)\cdot \hat{\mu}^{(d)}(Z) \right)^2.\label{eq:sigmasimple}
	\end{align}
    Now we  show the strong consistency of $\hat{\sigma}^2(nt,n\eta)$, i.e., $\hat{\sigma}^{2}(nt,n\eta)-\sigma^2 \xrightarrow[n\rightarrow\infty]{}0$ almost surely, as stated in Proposition \textcolor{black}{3.1}, returning to the time continuous setting but making use of formula \eqref{eq:sigmasimple}. We start by rewriting $\hat\sigma^2$ into a sum of six summands for each dimension $d$.
    \begin{align}
        \lefteqn{\hat{\sigma}^{2}(nt,n\eta)}\nonumber\\
        &=\frac{\sigma^2}{2n\eta-4}\sum_{j=1}^{n\eta}\sum_{d=1,2}\left(Z^{(d)}_{\lfloor nt\rfloor+j}-\bar{Z}^{(d)}-\frac{1}{2}w_{\mu}(n\eta,j)\cdot \hat{\mu}^{(d)}(Z) \right)^2\nonumber\\
		&=\sum_{d=1,2}\frac{\sigma^2}{2n\eta-4}\sum_{j=1}^{n\eta}\Biggl(Z^{(d)^2}_{\lfloor nt\rfloor+j}+\bar{Z}^{(d)^2}+\left(\frac{1}{2}w_{\mu}(n\eta,j)\cdot \hat{\mu}^{(d)}(Z)\right)^2 \nonumber\\
		&\qquad~+2Z^{(d)}_{\lfloor nt\rfloor+j}\left(-\bar{Z}^{(d)}\right) +2Z^{(d)}_{\lfloor nt\rfloor+j}\left(-\frac{1}{2}w_{\mu}(n\eta,j)\cdot \hat{\mu}^{(d)}(Z)\right)\nonumber\\
        &\qquad~+2\left(-\bar{Z}^{(d)}\right)\left(-\frac{1}{2}w_{\mu}(n\eta,j)\cdot \hat{\mu}^{(d)}(Z)\right)  \Biggr)\nonumber\\
		&=\sum_{d=1,2}\frac{\sigma^2}{2n\eta-4}\sum_{j=1}^{n\eta}\biggl(Z^{(d)^2}_{\lfloor nt\rfloor+j}+\bar{Z}^{(d)^2}+\frac{1}{4}w_{\mu}(n\eta,j)^2\hat{\mu}^{(d)}(Z)^2\nonumber\\
        &\qquad~-2Z^{(d)}_{\lfloor nt\rfloor+j}\bar{Z}^{(d)} -w_{\mu}(n\eta,j)Z^{(d)}_{\lfloor nt\rfloor+j}\hat{\mu}^{(d)}(Z) +w_{\mu}(n\eta,j)\bar{Z}^{(d)}\hat{\mu}^{(d)}(Z)  \biggr)\nonumber\\
        &=\sum_{d=1,2}\Bigg(\underbrace{\sum_{j=1}^{n\eta}\frac{\sigma^2}{2n\eta-4}Z^{(d)^2}_{\lfloor nt\rfloor+j}}_{=:S_1} +\underbrace{\sum_{j=1}^{n\eta}\frac{\sigma^2}{2n\eta-4}\bar{Z}^{(d)^2}}_{=:S_2}\nonumber\\
        &\qquad~+ \underbrace{\sum_{j=1}^{n\eta}\frac{\sigma^2}{2n\eta-4}\frac{1}{4}w_{\mu}(n\eta,j)^2\hat{\mu}^{(d)}(Z)^2}_{=:S_3}- \underbrace{\sum_{j=1}^{n\eta}\frac{\sigma^2}{2n\eta-4}2Z^{(d)}_{\lfloor nt\rfloor+j}\bar{Z}^{(d)}}_{=:S4}\nonumber\\
		&\qquad~ - \underbrace{\sum_{j=1}^{n\eta}\frac{\sigma^2}{2n\eta-4}w_{\mu}(n\eta,j)Z^{(d)}_{\lfloor nt\rfloor+j}\hat{\mu}^{(d)}(Z)}_{=:S_5}\nonumber\\
        &\qquad~+\underbrace{\sum_{j=1}^{n\eta}\frac{\sigma^2}{2n\eta-4}w_{\mu}(n\eta,j)\bar{Z}^{(d)}\hat{\mu}^{(d)}(Z)}_{=:S_6} \Bigg)\label{eq:hatsigmasquared_terms}
	\end{align}
	Since the summands in the two dimensions $d=1,2$ are independent, we focus on one dimension. As $\hat{\sigma}^{2}(nt,n\eta)$ is unbiased we can split the expectation $\sigma^2$ into the expectations of each of the summands. We show that each summand in the brackets in \eqref{eq:hatsigmasquared_terms} centred with its expectation converges to $0$ almost surely.\\
	For summand 1 we have
	\begin{align*}
		S_1-\mathbb{E}[S_1] &= \sum_{j=1}^{n\eta}\frac{\sigma^2}{2n\eta-4}\left(Z^{(d)^2}_{\lfloor nt\rfloor+j} - \mathbb{E}\bigg[Z^{(d)^2}_{\lfloor nt\rfloor+j}\bigg]\right)\\
        &= \sum_{j=1}^{n\eta}v(n,j)\left(\underbrace{Z^{(d)^2}_{\lfloor nt\rfloor+j} - \mathbb{E}\bigg[Z^{(d)^2}_{\lfloor nt\rfloor+j}\bigg]}_{\mathbb{E}[*]=0}\right),
	\end{align*}
	where we here define $v(n,j):=\frac{\sigma^2}{2n\eta-4}$. According to Proposition \textcolor{black}{ B.1}, we need to verify $\sum_{j=1}^{m_n}v(n,j)^2 {=} \mathrm{O}(n^{-1})$. This holds because
	\begin{align*}
		\sum_{j=1}^{n\eta}v(n,j)^2 = \sum_{j=1}^{n\eta}\left(\frac{\sigma^2}{2n\eta-4}\right)^2 = \sigma^4\frac{n\eta}{(2n\eta-4)^2} \kgs{=} \mathrm{O}(n^{-1}).
	\end{align*}
	For summand 2, we use $\mathbb E \bar{Z}^{(d)^2} = 1/(n\eta)$ and get
	\begin{align*}
		S_2-{\mathbb{E}[S_2]} 
		&= \frac{n\eta\sigma^2}{2n\eta-4} \left(\bar{Z}^{(d)^2}-\mathbb E \bar{Z}^{(d)^2}\right)\\
        &= \frac{n\eta \sigma^2}{2n\eta-4} \bar{Z}^{(d)^2}-\frac{\sigma^2}{2n\eta-4}\xrightarrow[n\rightarrow\infty]{}0 \text{ a.s.}
	\end{align*}
    For summand 3 we first rewrite $S_3$ as
	\begin{align*}
		S_3 
		&= \sum_{j=1}^{n\eta}\frac{\sigma^2}{2n\eta-4}\frac{1}{4}w_{\mu}(n\eta,j)^2\left(\frac{6}{(n\eta)^3-n\eta}\sum_{k=1}^{n\eta}w_{\mu}(n\eta,k)Z^{(d)}_{\lfloor nt\rfloor+k}\right)^2\\
		&= \sum_{k,l=1}^{n\eta}\underbrace{v(n,k,l)}_{\text{ in of Proposition \ref{proposition:SLLNnew_supp}}}Z^{(d)}_{\lfloor nt\rfloor+k}Z^{(d)}_{\lfloor nt\rfloor+l},
	\end{align*}
    with 
    \begin{align*}
        v(n,k,l) := \frac{6^2\sigma^2}{4(2n\eta-4)((n\eta)^3-n\eta)^2}\left(\sum_{j=1}^{n\eta}w_{\mu}(n\eta,j)^2\right) w_{\mu}(n\eta,k)w_{\mu}(n\eta,l)
    \end{align*}
    such that we have
    \begin{align*}
        S_3-\mathbb E [S_3] &= \sum_{k,l=1}^{n\eta}v(n,k,l)\left(Z^{(d)}_{\lfloor nt\rfloor+k}Z^{(d)}_{\lfloor nt\rfloor+l}-\mathbb E\left[Z^{(d)}_{\lfloor nt\rfloor+k}Z^{(d)}_{\lfloor nt\rfloor+l}\right]\right).
    \end{align*}
    To show almost sure convergence to $0$ we need to verify that the term in \eqref{eq:langerterm} is in $\mathrm{O}(n^{-2})$. We thus note that for the present definition of $v(n,k,l)$, the term in \eqref{eq:langerterm} corresponds to the following term
	\begin{align}
		  \lefteqn{\underbrace{\frac{6^8\sigma^8}{4^4(2n\eta-4)^4((n\eta)^3-n\eta)^8}\left(\frac{(n\eta)^3-n\eta)}{3}\right)^4}_{{=} \mathrm{O}(n^{-28}n^{12})=\mathrm{O}(n^{-16})}} \nonumber\\
		& ~\times \bigg[\underbrace{\sum_{j,k=1}^{n\eta}w_{\mu}(n\eta,j)^4w_{\mu}(n\eta,k)^4}_{(n\eta)^2\cdot \mathrm{O}(n^8)=\mathrm{O}(n^{10})}\nonumber\\
        & \qquad ~+\underbrace{\sum_{\substack{i,j,k,l=1\\ (i,j)\neq(k,l)}}^{n\eta}w_{\mu}(n\eta,i)^2w_{\mu}(n\eta,j)^2w_{\mu}(n\eta,k)^2w_{\mu}(n\eta,l)^2}_{{=}(n\eta)^4\cdot \mathrm{O}(n^8)=\mathrm{O}(n^{12})}\nonumber\\
		& \qquad ~+ \underbrace{\sum_{\substack{i,j,k,l=1\\ (i,j)\neq(k,l)}}^{n\eta}w_{\mu}(n\eta,i)^3w_{\mu}(n\eta,j)^3w_{\mu}(n\eta,k)w_{\mu}(n\eta,l)}_{{=} (n\eta)^4\cdot \mathrm{O}(n^8)=\mathrm{O}(n^{12})}\nonumber\\
		& \qquad ~+  \underbrace{\sum_{\substack{i,j,k,l,o,p=1\\ (i,j)\neq(k,l)\\ (i,j)\neq(o,p)\\ (o,p)\neq(k,l)}}^{n\eta} w_{\mu}(n\eta,i)w_{\mu}(n\eta,j)w_{\mu}(n\eta,k)w_{\mu}(n\eta,l)w_{\mu}(n\eta,o)^2w_{\mu}(n\eta,p)^2}_{{=}(n\eta)^6\cdot \mathrm{O}(n^8)=\mathrm{O}(n^{14})}\nonumber\\
		& \qquad ~+  \sum_{\substack{i,j,k,l,o,p,r,s=1;\\ (i,j),(k,l),\\(o,p),(r,s)\\ \text{ pairwise different}}}^{n\eta} \prod_{q=i,j,k,l,o,p,r,s}w_{\mu}(n\eta,q)\bigg]\label{eq:langerterm_summand3}.
\end{align}
All sums in the latter display except for the last one are easily upper bounded and seen to be in $\mathrm{O}(n^{14}/n^{16})=\mathrm{O}\left(n^{-2}\right)$. The last sum requires an exact, tedious calculation. We used the Python module Sympy and obtained  $\mathrm{O}(n^{12}/n^{16})$.

For summand 4 we use  $\mathbb{E}[Z^{(d)}_{\lfloor nt\rfloor+i}\bar{Z}] = 1/(n\eta)$ and get
	\begin{align}
		S_4 - \mathbb{E}[S_4] 
		&= \sum_{j=1}^{n\eta}\frac{\sigma^2}{n\eta-2}\left(Z^{(d)}_{\lfloor nt\rfloor+j}\bar{Z}^{(d)} - \frac{1}{n\eta}\right) \xrightarrow[n\rightarrow\infty]{}0\label{eq:consistencyhatvariance_summand4}
	\end{align}
by the conventional SLLN. 
 
For summand 5 we again first rewrite $S_5$ as
	\begin{align*}
	S_5 
		&= \sum_{i=1}^{n\eta}\frac{\sigma^2}{2n\eta-4}w_{\mu}(n\eta,i)Z^{(d)}_{\lfloor nt\rfloor+i}\left(\frac{6}{(n\eta)^3-n\eta}\sum_{j=1}^{n\eta}w_{\mu}( n\eta,j)Z^{(d)}_{\lfloor nt\rfloor+j}\right),
  \end{align*}
  which yields
\begin{align*}
		S_5-\mathbb E[S_5]&= \sum_{i,j=1}^{n\eta} v(n,i,j)\left(Z^{(d)}_{\lfloor nt\rfloor+i}Z^{(d)}_{\lfloor nt\rfloor+j}- \mathbb E \left[Z^{(d)}_{\lfloor nt\rfloor+i}Z^{(d)}_{\lfloor nt\rfloor+j}\right]\right)
	\end{align*}
	where we define 
    \begin{align*}
        v(n,i,j):=\frac{6\sigma^2}{(2n\eta-4)((n\eta)^3-n\eta)}w_{\mu}(n\eta,i)w_{\mu}(n\eta,j).
    \end{align*} 
    Now we are in the setting of the SLLN as in Proposition \ref{proposition:SLLNnew_supp} because the summands are not necessarily independent. Therefore we need to verify that the term in \eqref{eq:langerterm} is in $\mathrm{O}(n^{-2})$. We define the expression in the square brackets in \eqref{eq:langerterm_summand3} as $L$, for which we showed $L\kgs{=} \mathrm{O}(n^{14})$ and find that the term in \eqref{eq:langerterm} in the current setting can be written as
	\begin{align*}
		\underbrace{\frac{6^4\sigma^8}{(2n\eta-4)^4((n\eta)^3-n\eta)^4}}_{{=} \mathrm{O}(n^{-16})}\cdot \underbrace{L}_{{=} \mathrm{O}(n^{14})}{=}\mathrm{O}(n^{-2}).
	\end{align*}
For summand 6 we first observe that 
\begin{align}
		\mathbb{E}\bigg[\bar{Z}^{(d)}\hat{\mu}^{(d)}(Z) \bigg] 
		&= \frac{6}{n\eta((n\eta)^3-n\eta)}\sum_{k=1}^{n\eta}w_{\mu}(n\eta,k)\underbrace{\mathbb{E}\bigg[Z^{(d)^2}_{\lfloor nt\rfloor+k}\bigg]}_{=1}\nonumber\\
        &= \frac{6}{n\eta((n\eta)^3-n\eta)}\underbrace{\sum_{k=1}^{n\eta}w_{\mu}(n\eta,k)}_{=0}=0, \label{eq:expectationnebenrechnung}
	\end{align}
such that	
\begin{align*}
		{S_6 - \mathbb{E}[S_6]}
		= \sum_{i=1}^{n\eta}\frac{\sigma^2}{2n\eta-4}w_{\mu}(n\eta,i)\bar{Z}^{(d)}\hat{\mu}^{(d)}(Z)
 	\end{align*}
This term vanishes because 
  \begin{align*}
\sum_{i=1}^{n\eta}w_{\mu}(n\eta,i) = 0.
  \end{align*}
	
In total, each summand of $\hat{\sigma}^2(nt,n\eta)-\sigma^2$ converges to $0$ almost surely, which implies
	\begin{align*}
		\hat{\sigma}^{2}(nt,n\eta) \xrightarrow[n\rightarrow\infty]{} \sigma^{2} \text{   a.s.}
	\end{align*}

\end{document}